\documentclass[oneside, 12pt]{amsart}
\usepackage{amssymb,amsmath,amsfonts,amscd,amsthm,amsxtra,latexsym,scalerel,url,hyperref, mathtools}
\usepackage[utf8]{inputenc} 
\usepackage[T1]{fontenc}    
\usepackage{textcase}
\usepackage{stmaryrd}
\theoremstyle{plain}
\newtheorem{theorem}{Theorem}
\newtheorem{proposition}{Proposition}
\newtheorem{lemma}{Lemma}
\newtheorem{corollary}{Corollary}
\theoremstyle{remark}
\newtheorem{definition}{Definition}
\newtheorem{remark}{Remark}
\newtheorem{example}{Example}

\newcommand{\GL}{\mathrm{GL}}

\newcommand{\bbS}{\mathbb{S}}

\newcommand{\ev}{\mathrm{ev}}

\newcommand{\Sym}{\mathrm{Sym}}

\newcommand{\g}[2]{\langle\,#1,#2\,\rangle}
\renewcommand{\gg}{\g{\,\cdot\,}{\,\cdot\,}}

\newcommand{\bbR}{\mathbb{R}}
\newcommand{\C}{\mathbb{C}}

\newcommand{\N}{\mathbb{N}}

\newcommand{\R}{\mathrm{R}}
\newcommand{\cyclic}{\mathrm{cyclic}}

\newcommand{\Id}{\mathrm{Id}}

\newcommand{\SL}{\mathrm{SL}}

\newcommand{\Kern}{\mathrm{Kern}}

\newcommand{\rmS}{\mathrm{S}}

\newcommand{\bbA}{\mathbb{A}}
\newcommand{\bbE}{\mathbb{E}}

\newcommand{\scrA}{\mathcal{A}}

\newcommand{\scrC}{\mathcal{C}}

\newcommand{\scrM}{\mathcal{M}}
\newcommand{\scrQ}{\mathcal{Q}}
\newcommand{\scrR}{\mathcal{R}}

\newcommand{\scrO}{\mathcal{O}}

\newcommand{\jet}{\mathrm{jet}}

\newcommand{\rmR}{\mathrm{R}}

\DeclareMathOperator{\Endop}{End}
\newcommand{\End}[1]{\Endop\!\left(#1\right)}
\newcommand{\Endp}[1]{\Endop_{+}\!\left(#1\right)}

\title{The Jet Isomorphism Theorem of Riemannian Geometry}
\author{Tillmann Jentsch}
\address{Unidad Cuernavaca del Instituto de Matem\'aticas, Universidad
Nacional Aut\'onoma de M\'exico, Avenida Universidad s/n, Lomas de Chamilpa,
62210 Cuernavaca, Morelos, MEXICO (Guest visitor)}
\email{tilljentsch@gmail.com}


\begin{document}\sloppy
\maketitle

\begin{abstract}
  A classical theorem of Riemannian geometry, due in its original
  form to Cartan, states that the Taylor expansion of the metric in
  geodesic normal coordinates is a universal formal power series
  involving only the symmetrizations of the iterated covariant
  derivatives of the curvature tensor; this is known as the jet
  isomorphism theorem. In particular, it is in principle possible
  to reconstruct the jet of the curvature tensor from its symmetrization
  in geodesic normal coordinates, although this would certainly
  result in an unwieldy computation. In this paper we achieve
  the same goal by coordinate--free calculations, using only the intrinsic definition of the relevant Young symmetrizers.
\end{abstract}

\section{Overview}

Let \(M\) be a smooth manifold equipped with a nondegenerate symmetric tensor
field \(g\) of type \((0,2)\). The pair \((M,g)\) is called a
semi-Riemannian or pseudo-Riemannian manifold. In what follows, we shall
simply refer to \((M,g)\) as a Riemannian manifold.

A Riemannian manifold \((M,g)\) admits a unique torsion-free connection
\(\nabla\) satisfying \(\nabla g = 0\); this is the Levi–Civita connection.
Associated with \(\nabla\) is the Riemann curvature tensor
\[
  \rmR(X,Y)Z \coloneqq \nabla_X \nabla_Y Z - \nabla_Y \nabla_X Z
  - \nabla_{[X,Y]} Z
\]
for all \(X,Y,Z \in \Gamma(TM)\), where \(\Gamma(TM)\) denotes the space of
smooth vector fields on \(M\). For \(k \ge 0\), we denote by \(\nabla^k \rmR\)
its \(k\)-fold iterated covariant derivative. For \(X,Y \in T_pM\) we put
\begin{equation}\label{eq:k-ter_Jacobi}
\scrR^k(X)Y \coloneqq \nabla^{k}_{X,\ldots,X}\,\rmR(Y,X)X
\end{equation}
called the \emph{symmetrized \(k\)th covariant derivative of the curvature
tensor}. By definition, \(\scrR^k|_p\) is a polynomial map
\[
\scrR^k|_p \colon T_pM \longrightarrow \Endp{T_pM}, \quad
X \longmapsto \scrR^k(X)
\]
of degree \(k+2\) on \(T_pM\) with values in \(\Endp{T_pM}\), the space of
symmetric endomorphisms of \(T_pM\).

In principle it is possible to reconstruct the \(k\)-jet
\[
\nabla^{\le k}|_p \rmR \coloneqq
\big(\rmR|_p,\, \nabla|_p \rmR,\, \ldots,\, \nabla^k|_p \rmR\big)
\]
of the curvature tensor from its symmetrization
\[
\scrR^{\le k}|_p \coloneqq
\big(\scrR^0|_p,\, \scrR^1|_p,\, \ldots,\, \scrR^k|_p\big)
\]
at an arbitrary point \(p\) via the following classical result.

\bigskip
\begin{theorem}[{\cite[Appendix~II]{ABP}}, \cite{Gray}]
\label{th:taylor_series}
Let \(p \in M\), let \(\gg \coloneqq g|_p\) denote the inner product on
\(T_pM\), and \(\exp^M_p\colon U \to M\) be the exponential map, defined on
an open star-shaped neighborhood \(U \subset T_pM\) of the origin. Let
\(\tilde g\) denote the pullback of \(g\) under \(\exp^M_p\), i.e.,~the metric
tensor in geodesic normal coordinates at \(p\). Then there exist universal
noncommutative polynomials \(Q_k\) of degree \(k\) in a countable set of free
variables such that
\begin{equation}\label{eq:Taylor1}
\langle Y, Z \rangle +
\sum_{j=2}^k \frac{1}{j!}\,\big\langle
Q_j\big(\scrR^0(X), \scrR^1(X), \ldots\big) Y, Z \big\rangle
\end{equation}
is the Taylor polynomial of order \(k\) for the function
\(X \mapsto \tilde g(X)_{Y,Z}\) on \(T_pM\) for all \(Y,Z \in T_pM\) and
\(k \ge 0\).
\end{theorem}

For example,
\begin{equation}\label{eq:Taylor_of_deg_5}
\begin{aligned}
&\langle Y, Z \rangle
 - \frac{1}{3}\,\langle \scrR^0(X)Y, Z \rangle
 - \frac{1}{6}\,\langle \scrR^1(X)Y, Z \rangle
 - \frac{1}{20}\,\langle \scrR^2(X)Y, Z \rangle
 + \frac{2}{45}\,\langle \scrR^0(X)\scrR^0(X)Y, Z \rangle \\
&\quad - \frac{1}{90}\,\langle \scrR^3(X)Y, Z \rangle
 + \frac{1}{45}\,\big\langle
 \scrR^0(X)\scrR^1(X)Y + \scrR^1(X)\scrR^0(X)Y, Z \big\rangle
\end{aligned}
\end{equation}
is the Taylor polynomial of order five of the metric tensor in geodesic
normal coordinates at any point \(p \in M\). In
Appendix~\ref{se:Taylor_expansion} we will explain in detail the notation
used in Theorem~\ref{th:taylor_series} and recall its proof. There we will
also give a recursive formula for the coefficients of the Taylor series of
the backward parallel transport map; see
Proposition~\ref{p:taylor_series_of_the_backward_parallel_transport}.

By the invariance of curvature jets under isometries and since \(\exp^M_p\)
is an anchored coordinate system based at \(p\), in the sense that
\(d(\exp^M_p)_0 = \Id_{T_pM}\), we obtain
\[
\bigl(\tilde{\nabla}^{\le k}\tilde{\rmR}\bigr)\big|_{0}
= \bigl(\nabla^{\le k}\rmR\bigr)\big|_{p}
\]
where \(\tilde{\nabla}^{\le k}\tilde{\rmR}\) denotes the \(k\)-jet of the
curvature tensor corresponding to the polynomial metric \(\tilde g\) on
\(T_pM\) defined by \eqref{eq:Taylor1}. Thus, in principle, one can recover
\((\nabla^{\le k}\rmR)|_p\) from its symmetrization by working in geodesic
normal coordinates. In practice, however, this would require knowledge of
the Levi–Civita connection \(\tilde{\nabla}\) of \(\tilde g\), the curvature
tensor \(\tilde{\rmR}\), and its iterated covariant derivatives, which does
not seem to yield a useful closed formula in any straightforward way.

One of the main goals of this exposition is therefore to find an explicit
recursive formula for \(\nabla^{\le k}|_p \rmR\) in terms of
\(\scrR^{\le k}|_p\), while completely bypassing the Taylor expansion of the
metric in geodesic normal coordinates; see \eqref{eq:YS2}. In fact, we
subsequently also obtain a practical formula for the \(k\)-jet of the
curvature tensor of the metric in geodesic normal coordinates. For this, it
only remains to solve the equations
\(h_{k+2} = Q_{k+2}(\scrR^{\le k}|_p)\) from \eqref{eq:Taylor1} for
\(\scrR^{\le k}|_p\), which is part of the classical jet isomorphism theorem
restated in Theorem~\ref{th:jet_isomorphism_theorem}.

\subsection{The inverse of the jet symmetrization map}

As usual, we also set
\[
\nabla^k_{X_5,\ldots,X_{k+4}} \rmR_{X_1,X_2,X_3,X_4}
\coloneqq
\big\langle \nabla^k_{X_{5},\ldots,X_{k+4}} \rmR_{X_1,X_2,X_3},\,X_4 \big\rangle
\]
which means that the \((1,k+3)\)-tensor \(\nabla^k \rmR\) can also be regarded,
in a natural way, as a tensor of type \((0,k+4)\) for all \(k \ge 0\).
In the same vein, \(\scrR^k\), defined in \eqref{eq:k-ter_Jacobi}, can be viewed
as a section of \(\Sym^{k+2}T^*M \otimes \Sym^{2}T^*M\) characterized by
\(\scrR^k_{X,\ldots,X;\,Y,Z} = \langle \scrR^k(X)Y, Z \rangle\).

Let \(\rmS^\star_{\scaleto{
 \begin{array}{|c|c|c|c|c|c|}\cline{1-5}
  1 & 3 & 5 & \cdots & k\!+\!4\\
 \cline{1-5}
  2 & 4 & \multicolumn{1}{c}{\;\;\;}\\
 \cline{1-2}
\end{array}
}{15pt}}\) denote the Young symmetrizer associated with the standard Young tableau
\begin{equation}\label{eq:Fiedler_alternativ}
 \begin{array}{|c|c|c|c|c|c|}
  \cline{1-5}
  1 & 3 & 5 & \cdots & k+4\\
  \cline{1-5}
  2 & 4 & \multicolumn{2}{c}{\;\;\;} \\
  \cline{1-2}
 \end{array}
\end{equation}
of shape \((k+2,2)\). By definition,
\begin{equation}\label{eq:YS1}
 \rmS^\star_{\scaleto{
 \begin{array}{|c|c|c|c|c|c|}\cline{1-5}
  1 & 3 & 5 & \cdots & k\!+\!4\\
  \cline{1-5}
  2 & 4 & \multicolumn{1}{c}{\;\;\;}\\
  \cline{1-2}
 \end{array}
 }{15pt}}
 \nabla^k_{X_5,\ldots,X_{k+4}}\,\rmR_{X_1,X_2,X_3,X_4}
 =
 -2(k+2)!\,
 \rmS_{\scaleto{
 \begin{array}{|c|}\cline{1-1}
  1\\\cline{1-1}
  2\\\cline{1-1}
 \end{array}
 }{15pt}}\,
 \rmS_{\scaleto{
 \begin{array}{|c|}\cline{1-1}
  3\\\cline{1-1}
  4\\\cline{1-1}
 \end{array}
 }{15pt}}\,
 \scrR^k_{X_1,X_3,X_5,\ldots,X_{k+4};\,X_2,X_4}
\end{equation}
Here,
\(\rmS_{\scaleto{
 \begin{array}{|c|}\cline{1-1}
  1\\\cline{1-1}
  2\\\cline{1-1}
 \end{array}
 }{15pt}}\) and
\(\rmS_{\scaleto{
 \begin{array}{|c|}\cline{1-1}
  3\\\cline{1-1}
  4\\\cline{1-1}
 \end{array}
 }{15pt}}\) are the antisymmetrizers in the pairs of variables \(\{X_1,X_2\}\)
and \(\{X_3,X_4\}\), respectively. Then,
\[
\rmS_{\scaleto{
 \begin{array}{|c|}\cline{1-1}
  1\\\cline{1-1}
  2\\\cline{1-1}
 \end{array}
 }{15pt}}\,
\rmS_{\scaleto{
 \begin{array}{|c|}\cline{1-1}
  3\\\cline{1-1}
  4\\\cline{1-1}
 \end{array}
 }{15pt}}
 = \owedge \otimes \Id
\]
where \(\owedge \colon \Sym^{2}V^* \otimes \Sym^2 V^* \to \Sym^2(\Lambda^2 V^*)\)
denotes the classical Kulkarni–Nomizu product in the variables
\(\{X_1,X_2,X_3,X_4\}\):
\begin{equation}\label{eq:Kulkarni_Nomizu}
 (h \owedge \tilde h)_{X_1,\ldots,X_{4}}
 \coloneqq
  h_{X_1,X_3} \tilde h_{X_2,X_4}
 - h_{X_2,X_3} \tilde h_{X_1,X_4}
 - h_{X_1,X_4} \tilde h_{X_2,X_3}
 + h_{X_2,X_4} \tilde h_{X_1,X_3}
\end{equation}
and \(\Id\) is the identity map on covariant \(k\)-tensors in the variables
\(X_5,\ldots,X_{k+4}\).

For the moment, assume that \(\rmR\) and its first \(k-1\) covariant derivatives
vanish at a given point: \(\nabla^{\ell}|_p \rmR = 0\) for \(0 \le \ell \le k-1\).
Then \(\nabla^{\le k}|_p \rmR\) is a \emph{linear} \(k\)-jet;
see Definition~\ref{de:algebraic_curvature_jets}\,(d). In this case, the Young
projection formula
\begin{equation}\label{eq:Young_Symmetrizer}
 \nabla^k_{X_5,\ldots,X_{k+4}}\,\rmR_{X_1,X_2,X_3,X_4}
 =
 \frac{1}{h_k}\,
 \rmS^\star_{\scaleto{
 \begin{array}{|c|c|c|c|c|c|}\cline{1-5}
  1 & 3 & 5 & \cdots & k\!+\!4\\
  \cline{1-5}
  2 & 4 & \multicolumn{1}{c}{\;\;\;}\\
  \cline{1-2}
 \end{array}
 }{15pt}}
 \nabla^k_{X_5,\ldots,X_{k+4}}\,\rmR_{X_1,X_2,X_3,X_4}
\end{equation}
holds, where \(h_k = 2\,k!\,(k+2)(k+3)\) is the product of the hook lengths of
the Young frame underlying \eqref{eq:Fiedler_alternativ}; see
Section~\ref{se:YS}. From \eqref{eq:YS1}--\eqref{eq:Young_Symmetrizer} we see
how a linear \(k\)-jet can be reconstructed from its symmetrization:
\begin{equation}\label{eq:Fiedler}
 \nabla^k \rmR = -\,\frac{k+1}{k+3}\,(\owedge \otimes \Id)\,\scrR^k
\end{equation}

Although this part of the theory is well established (indeed, it is also one
of the main arguments in the proof of the classical jet isomorphism theorem),
in Section~\ref{se:Fiedler} we nevertheless present an elementary proof of
\eqref{eq:Young_Symmetrizer}, using only the direct definition of
\(\rmS^\star_{\scaleto{
 \begin{array}{|c|c|c|c|c|c|}\cline{1-5}
  1 & 3 & 5 & \cdots & k\!+\!4\\
  \cline{1-5}
  2 & 4 & \multicolumn{1}{c}{\;\;\;}\\
  \cline{1-2}
 \end{array}
 }{15pt}}\) (as explained earlier) together with the two Bianchi identities.
In other words, we show by direct calculation that Weyl’s construction of an
irreducible representation of \(\SL(n,\C)\) of highest weight \((k+2,2)\)
contains the intersection of the two Bianchi identities, which is the
nontrivial part of Theorem~\ref{th:char1} in this special case.\footnote{It
would be interesting to know whether Theorem~\ref{th:char1} can be proved in a
similar way for arbitrary Young frames.}

To understand the general case (i.e.,~when \(\nabla^{\le k}|_p \rmR\) is not
necessarily linear), in Section~\ref{se:general_case} we consider the
symmetrized iterated covariant derivative
\[
\jet^k_{X_1,\ldots,X_k} \rmR
\coloneqq
\frac{1}{k!} \sum_{\sigma \in \rmS_k}
\nabla^k_{X_{\sigma(1)},\ldots,X_{\sigma(k)}} \rmR
\]
where \(\rmS_k\) denotes the symmetric group. By the Ricci identity,
\(\jet^k|_p \rmR = \nabla^k|_p \rmR\) holds for every linear \(k\)-jet.
Therefore, we may rewrite \eqref{eq:Young_Symmetrizer} as
\begin{equation}\label{eq:Young_Symmetrizer_2}
\Big(
 \frac{1}{k!}\,
 \rmS^\star_{\scaleto{
 \begin{array}{|c|c|c|c|c|c|}\cline{1-5}
  1 & 3 & 5 & \cdots & k\!+\!4\\
  \cline{1-5}
  2 & 4 & \multicolumn{1}{c}{\;\;\;}\\
  \cline{1-2}
 \end{array}
}{15pt}}
\nabla^k_{X_5,\ldots,X_{k+4}}
 - 2(k+2)(k+3)\,\jet^k_{X_5,\ldots,X_{k+4}}
\Big)
 \rmR_{X_1,X_2,X_3,X_4}
 = 0
\end{equation}

Our detailed approach to \eqref{eq:Young_Symmetrizer} in
Section~\ref{se:Fiedler} allows us to determine precisely how
\eqref{eq:Young_Symmetrizer_2} must be modified when the hypothesis
\(\nabla^\ell|_p \rmR = 0\) for \(\ell < k\) is dropped;
see Proposition~\ref{p:Young_Symmetrizer_modifiziert_1}. In this case,
the left-hand side of \eqref{eq:Young_Symmetrizer_2} does not necessarily
vanish but is given by the following term:
\begin{equation}\label{eq:Young_Symmetrizer_modifiziert_1}
\begin{aligned}
&\sum_{A=1}^k
 \rmS^\star_{\scaleto{
 \begin{array}{|c|c|c|}\cline{1-3}
  1 & 3 & A\!+\!4\\
  \cline{1-3}
  2 & 4 & \multicolumn{1}{c}{\;\;\;}\\
  \cline{1-2}
 \end{array}
}{15pt}}
\big(
 \jet^k_{X_5,\ldots,X_{k+4}}
 - \jet^{k-1}_{X_5,\ldots,\hat X_{A+4},\ldots,X_{k+4}}\,
   \nabla_{X_{A+4}}
\big) \rmR_{X_1,X_2,X_3,X_4} \\
&\quad+ \sum_{\substack{A,B=1 \\ A < B}}^k
 \rmS^\star_{\scaleto{
 \begin{array}{|c|c|}\cline{1-2}
  1 & 3\\
  \cline{1-2}
  2 & 4\\
  \cline{1-2}
 \end{array}
}{15pt}}
\big(
 \jet^k_{X_1,X_3,X_5,\ldots,\hat X_{A+4},\ldots,\hat X_{B+4},\ldots,X_{k+4}}
 - \jet^{k-2}_{X_5,\ldots,\hat X_{A+4},\ldots,\hat X_{B+4},\ldots,X_{k+4}}\,
   \nabla^2_{X_1,X_3}
\big) \rmR_{X_{A+4},X_2,X_{B+4},X_4} \\
&\quad+ \sum_{\substack{A,B=1 \\ A \neq B}}^k
\big(
 2\,
 \rmS_{\scaleto{
 \begin{array}{|c|}\cline{1-1}
  1\\\cline{1-1}
  2\\\cline{1-1}
 \end{array}
}{15pt}}
 \jet^{k-2}_{X_5,\ldots,\hat X_{A+4},\ldots,\hat X_{B+4},\ldots,X_{k+4}}\,
 \rmR_{X_1,X_{A+4}}\,\rmR_{X_{B+4},X_2,X_3,X_4} \\
&\hspace{5.6em}
 -
 \rmS_{\scaleto{
 \begin{array}{|c|}\cline{1-1}
  1\\\cline{1-1}
  2\\\cline{1-1}
 \end{array}
}{15pt}}
 \rmS_{\scaleto{
 \begin{array}{|c|}\cline{1-1}
  3\\\cline{1-1}
  4\\\cline{1-1}
 \end{array}
}{15pt}}
 \jet^{k-2}_{X_5,\ldots,\hat X_{A+4},\ldots,\hat X_{B+4},\ldots,X_{k+4}}\,
 \rmR_{X_1,X_3}\,\rmR_{X_{A+4},X_2,X_{B+4},X_4}
\big)
\end{aligned}
\end{equation}

Here, differential operators act on sections of induced vector bundles,
i.e.,~the Leibniz rule is not yet incorporated. It remains to understand
terms of the form \((\jet^k - \jet^\ell \nabla^{k-\ell}) \rmR\) for
\(0 \le \ell \le k\), which is a less specific problem. Its general
solution is given in Proposition~\ref{p:Simple_Jet-Formula} of
Section~\ref{se:jet_theory}. For this, \(\psi \coloneqq \rmR\) could in fact
be any section of some vector bundle \(\bbE\) with a linear connection
\(\nabla^\bbE\).

Explicit calculations for \(k \le 5\) are provided in
Appendix~\ref{se:special_jet_formula}. We ultimately find that the correct
modification of \eqref{eq:Fiedler} is
\begin{equation}\label{eq:YS2}
 \nabla^k \rmR + \frac{k+1}{k+3}(\owedge \otimes \Id)\,\scrR^k
 = B(\nabla^{\le k-2} \rmR)
\end{equation}
where \(B(\nabla^{\le k-2} \rmR)\) is a covariant \((k+4)\)-tensor that is a
quadratic expression in the \((k-2)\)-jet, which can be determined
explicitly; cf.~Corollary~\ref{co:B}. For small values of \(k\), see
Example~\ref{ex:tilde_B} together with Example~\ref{ex:k=2345} in
Section~\ref{se:explicit_calculations}.

\subsection{Outlook: Natural equations for jets of the curvature tensor}

Looking at the structurally involved algebraic properties that distinguish
\(\nabla^{\le k} \rmR\) (summarized in Definition~\ref{de:algebraic_curvature_jets}
below), it seems more advantageous to work instead with the symmetrized jet
\(\scrR^{\le k}\). In this formulation, the two Bianchi identities become a single
equation \eqref{eq:formal_Jacobi_operator}, meaning that \(\scrR^{\le k}\) is a
section of the graded vector bundle of algebraic symmetrized jets
\[
\scrC_\bullet(M) \coloneqq \bigoplus_{j=0}^\infty \scrC_j(TM)
\]
Clearly, the fiber \(\scrC_\bullet(M)_p = \bigoplus_{j=0}^\infty \scrC_j(T_pM)\) is
not only a vector space but also a graded module over the polynomial ring
\(\Sym^\bullet T^*_pM \coloneqq \bigoplus_{j=0}^\infty \Sym^j T^*_pM\). In particular,
the symmetrized \(k\)-jet is amenable to methods of commutative algebra.

In this setting, it is natural to study equations of the type
\begin{equation}\label{eq:JR}
 \scrR^k = -\sum_{i=1}^k a_i\,\scrR^{k-i}
\end{equation}
between polynomial functions of degree \(k+2\) on \(T_pM\) with values in
\(\Endp{T_pM}\), where each \(a_i\) is a polynomial function of degree \(i\) on
\(T_pM\). By definition, \eqref{eq:JR} means that
\begin{equation}\label{eq:JR_2}
 \scrR^k(X) = -\sum_{i=1}^k a_i(X)\,\scrR^{k-i}(X)
\end{equation}
for all \(X \in T_pM\). In \cite{JW2} we called \eqref{eq:JR} a \emph{Jacobi relation};
however, this terminology may not be optimal. Therefore, in \cite{J3} and \cite{J4},
a normed polynomial
\[
P(\lambda) \coloneqq \lambda^k + \sum_{i=1}^k a_i\,\lambda^{k-i}
\]
with coefficients \(a_i\) in the vector space of polynomial functions of degree \(i\) on
\(T_pM\) is called \emph{admissible} if \eqref{eq:JR} holds with the same \(k\) and the same
coefficients \(a_i\). For instance, admissible polynomials exist pointwise by Hilbert’s basis
theorem. Furthermore, for every compact real analytic Riemannian space, there exists a globally
admissible polynomial; that is, there are smooth sections \(a_i\) of the vector bundles of
polynomial functions of degree \(i\) on the various tangent spaces such that \eqref{eq:JR} holds
at each point of \(M\); cf.~the appendix of \cite{J3}.

As an application of the jet isomorphism theorem, \eqref{eq:JR} holds if and only if the curvature
tensor satisfies the explicit partial differential equation
\begin{equation}\label{eq:PDE}
 \nabla^k \rmR = A_{a_1,\ldots,a_k}(\nabla^{\le k-1} \rmR)
 + B(\nabla^{\le k-2} \rmR)
\end{equation}
where \(B(\nabla^{\le k-2} \rmR)\) is the term from the right-hand side of
\eqref{eq:YS2}, and \(A_{a_1,\ldots,a_k}(\nabla^{\le k-1} \rmR)\) is the section of the
bundle \(\scrC^\star_k(TM)\) of linear curvature \(k\)-jets over \(M\) defined by
\begin{equation}\label{eq:A}
\begin{aligned}
 &A_{a_1,\ldots,a_k}(\nabla^{\le k-1} \rmR)_{X_1,\ldots,X_{k+4}} \\
 &\quad\coloneqq
 -\frac{1}{h_k}\,
 \rmS^\star_{\scaleto{
 \begin{array}{|c|c|c|c|c|c|}\cline{1-5}
  1 & 3 & 5 & \cdots & k\!+\!4\\
  \cline{1-5}
  2 & 4 & \multicolumn{1}{c}{\;\;\;}\\
  \cline{1-2}
 \end{array}
 }{15pt}}
 \sum_{i=1}^k a_i(X_{k-i+5},\ldots,X_{k+4})\,
 \nabla^{k-i}_{X_5,\ldots,X_{k-i+4}} \rmR_{X_1,X_2,X_3,X_4}
\end{aligned}
\end{equation}
The proof of \eqref{eq:A} is straightforward using \eqref{eq:YS1} and \eqref{eq:YS2};
the details are left to the reader.

To take a broader view, recall that perhaps the most natural higher-order PDE one might
imagine for the curvature tensor, namely \(\nabla^k \rmR = 0\), already implies
\(\nabla \rmR = 0\)—that is, the manifold is Riemannian symmetric—whenever \(M\) is
complete; cf.~\cite{NO}. Since compact real analytic Riemannian spaces occur in abundance
\cite{MO}, one may view \eqref{eq:PDE} as a substitute for \(\nabla^k \rmR = 0\) that
still incorporates a rich variety of interesting examples.

\section{Revising the well-known parts of the Jet Isomorphism Theorem}

Let \(\bbE\) be a vector bundle over \(M\) equipped with a connection
\(\nabla\) (e.g., a tensor bundle with the connection induced by the
Levi–Civita connection). Following \cite[p.~23]{W1}, define the higher
covariant derivatives \(\nabla^k\psi\) of a section
\(\psi\in\Gamma(\bbE)\) iteratively, for \(k\ge 0\), by
\[
  \nabla^{k+1}_{Y,X_1,\ldots,X_k}\psi
  := \nabla_Y \nabla^{k}_{X_1,\ldots,X_k}\psi
     - \sum_{i=1}^k \nabla^{k}_{X_1,\ldots,\nabla_Y X_i,\ldots,X_k}\psi
\]
Hence the \(k\)-jet \(\nabla^{\le k}\psi :=
(\psi,\nabla\psi,\ldots,\nabla^k\psi)\) is a section of
\(\bigoplus_{i=0}^k \bigotimes^i T^*M \otimes \bbE\). Since the
Levi–Civita connection is torsion-free, the Ricci identity
\[
  \nabla^2_{X,Y} - \nabla^2_{Y,X} = \rmR^\bbE_{X,Y}
\]
holds, where the curvature endomorphism \(\rmR^\bbE_{X,Y}\colon
\bbE_p\to\bbE_p\) acts by \(\psi \mapsto \rmR^\bbE_{X,Y}\psi\).
Therefore,
\begin{equation}\label{eq:Ricci_psi}
  \nabla^{k+\ell+2}_{X_1,\ldots,X_k,A,B,Y_1,\ldots,Y_\ell}\psi
  - \nabla^{k+\ell+2}_{X_1,\ldots,X_k,B,A,Y_1,\ldots,Y_\ell}\psi
  = \nabla^k_{X_1,\ldots,X_k}\rmR_{A,B}\nabla^{\ell}_{Y_1,\ldots,Y_\ell}\psi
\end{equation}
for all \(k,\ell\ge 0\). Here \(\rmR\) and \(\nabla^k\) denote,
respectively, the curvature tensor and the \(k\)-fold covariant
derivative with respect to the induced connections on
\(\otimes^\ell T^*M\otimes\bbE\) and
\(\otimes^\ell T^*M\otimes\bbE\otimes\Lambda^2 T^*M\).

The statement \eqref{eq:Ricci_psi} omits the Leibniz rule.
Incorporating it yields (cf.~\cite[(3.1), p.~23]{W1})
\begin{equation}\label{eq:Ricci_refined}
\begin{aligned}
  &\nabla^{k+\ell+2}_{X_1,\ldots,X_k,A,B,Y_1,\ldots,Y_\ell}\psi
   - \nabla^{k+\ell+2}_{X_1,\ldots,X_k,B,A,Y_1,\ldots,Y_\ell}\psi \\
  &\qquad= \sum_{r=0}^k \sum_{1\le \mu_1<\cdots<\mu_r\le k}
     \Big(
       (\nabla^r_{X_{\mu_1},\ldots,X_{\mu_r}}\rmR^\bbE_{A,B})
       \nabla^{k+\ell-r}_{X_1,\ldots,\hat X_{\mu_1},\ldots,\hat X_{\mu_r},\ldots,X_k,\,
                           Y_1,\ldots,Y_\ell}\psi \\
  &\qquad\qquad\qquad\qquad
       - \sum_{\nu=1}^{\ell}
         \nabla^{k+\ell-r}_{X_1,\ldots,\hat X_{\mu_1},\ldots,\hat X_{\mu_r},\ldots,X_k,\,
                            Y_1,\ldots,(\nabla^r_{X_{\mu_1},\ldots,X_{\mu_r}}\rmR)_{A,B}Y_\nu,\ldots,Y_\ell}
         \psi
     \Big)
\end{aligned}
\end{equation}
where hats indicate omitted arguments.

More specifically, the \(k\)-jet \(\nabla^{\le k}\rmR :=
(\rmR,\nabla\rmR,\ldots,\nabla^k\rmR)\) is a section of
\(\bigoplus_{i=0}^k \bigotimes^{i+4} T^*M\) and satisfies the
classical algebraic constraints (cf.~\cite[Def.~8.2]{FG}):

\bigskip
\begin{definition}\label{de:algebraic_curvature_jets}
\begin{enumerate}
\item We have \(\rmR|_p \in \Sym^2(\Lambda^2 V^*)\) for \(V:=T_pM\):
\[
  \rmR(X_1,X_2,Y_1,Y_2) = -\,\rmR(X_2,X_1,Y_1,Y_2) =
  \rmR(Y_1,Y_2,X_1,X_2)
\]
and the first Bianchi identity
\[
  \rmR(X_1,X_2,X_3,Y) + \rmR(X_2,X_3,X_1,Y) +
  \rmR(X_3,X_1,X_2,Y) = 0
\]
for all \(X_1,X_2,X_3,Y \in V\). Conversely, for any Euclidean vector
space \((V,\gg)\) these equations define the subspace
\(\scrC_0^\star(V)\subset \otimes^4 V^*\) of algebraic curvature
tensors (associated to some metric \(g\) on \(V\) with \(g_0=\gg\)).

\item For each \(X\in T_pM\), the \(4\)-tensor
\(X\,\lrcorner\,\nabla\rmR := \nabla_X\rmR\) is an algebraic curvature
tensor, and the second Bianchi identity holds:
\[
  \nabla_{X_1}\rmR(X_2,X_3,Y_1,Y_2)
  + \nabla_{X_2}\rmR(X_3,X_1,Y_1,Y_2)
  + \nabla_{X_3}\rmR(X_1,X_2,Y_1,Y_2) = 0
\]
for all \(X_1,X_2,X_3,Y_1,Y_2 \in V\). Thus, for any Euclidean
\((V,\gg)\) these properties define the space \(\scrC_1^\star(V)\) of
algebraic covariant derivatives of the curvature tensor (associated to
some metric \(g\) on \(V\) with \(g_0=\gg\)).

\item For \(k\ge 2\), \(\nabla^{\le k}|_p\rmR\) has the following
characteristic properties: for all \(1\le \ell \le k-1\) and
\(X_1,\ldots,X_\ell \in V\),
\[
  (X_1\otimes\cdots\otimes X_\ell)\,\lrcorner\,\nabla^{\ell+1}\rmR
  := \nabla^\ell_{X_1,\ldots,X_\ell}\nabla\rmR \in \scrC_1^\star(V)
\]
and, by the Ricci identity \eqref{eq:Ricci_psi},
\begin{equation}\label{eq:Ricci_R}
  \begin{aligned}
  &\nabla^\ell_{X_1,\ldots,X_{\ell_1},A,B,Y_1,\ldots,Y_{\ell_2}}\rmR
   - \nabla^\ell_{X_1,\ldots,X_{\ell_1},B,A,Y_1,\ldots,Y_{\ell_2}}\rmR
   \\
  &\qquad= \nabla^{\ell_1}_{X_1,\ldots,X_{\ell_1}}\rmR_{A,B}
            \nabla^{\ell_2}_{Y_1,\ldots,Y_{\ell_2}}\rmR
  \end{aligned}
\end{equation}
for all \(2\le \ell\le k\) and \(\ell_1,\ell_2\ge 0\) with
\(\ell=\ell_1+\ell_2+2\). By \eqref{eq:Ricci_refined}, this is an
intrinsic tensorial property of \(\nabla^{\le k}|_p\rmR\).
Conversely, given a Euclidean \((V,\gg)\), any element of
\(\bigoplus_{\ell=0}^k \bigotimes^{\ell+4} V^*\) with these properties
is called an \emph{algebraic \(k\)-jet of the curvature tensor}
(associated to some metric \(g\) on \(V\) with \(g_0=\gg\)).

\item Suppose \(\nabla^\ell|_p\rmR = 0\) for \(0\le \ell\le k-1\). We
call \(\nabla^{\le k}|_p\rmR\) \emph{linear}. By the Ricci identity,
\begin{equation}\label{eq:scrC_k_star}
  \nabla^k|_p\rmR \in
  \scrC_k^\star(V) := \Sym^k(V^*)\otimes \scrC_0^\star(V)\ \cap\
  \Sym^{k-1}(V^*)\otimes \scrC_1^\star(V)
\end{equation}
with \(V:=T_pM\). Conversely, for any Euclidean \((V,\gg)\) the
elements of \(\scrC_k^\star(V)\) are called \emph{linear algebraic
\(k\)-jets}.
\end{enumerate}
\end{definition}

By the jet isomorphism theorem proved below, every algebraic \(k\)-jet
is the actual jet \(\nabla^{\le k}|_0\rmR\) of the curvature tensor of
some Riemannian metric \(g\) on \(V\) with \(g|_0=\gg\). In particular,
the theorem applies to any algebraic curvature tensor \(\rmR|_0\), to
any algebraic covariant derivative \(\nabla|_0\rmR\), and more
generally to any linear algebraic \(k\)-jet whose only nonzero
component is \(\nabla^k|_0\rmR\).

\subsection{Geodesic normal coordinates}

In geodesic normal coordinates \(\exp^M_p : T_pM \to M\), the
geodesics of \((M,g)\) emanating from \(p\) become the straight lines
emanating from the origin of \(T_pM\). This suggests the following
definition:

\bigskip
\begin{definition}\label{de:geod_norm_coord}
 Let \(V\) be a vector space and \(U\) a star-shaped open
 neighborhood of the origin. A Riemannian metric \(\tilde{g} : U \to
 \Sym^2_{\mathrm{reg}} V^*\) is given in geodesic normal coordinates
 if the straight lines \(t \mapsto tX\) are geodesics for all \(X \in
 U\) and \(t \in [0,1]\).
\end{definition}

Here, \(\Sym^2_{\mathrm{reg}} V^*\) denotes the set of nondegenerate
symmetric bilinear forms on \(V^*\). An anchored coordinate system
\(f : T_pM \to M\) (i.e.,~such that \(d_p f = \Id\)) is the geodesic
normal coordinate system for a given metric \(g\) on \(M\) if and only
if \(\tilde{g} := f^* g\) is given in normal coordinates. The
following classical result tells us how to detect geodesic normal
coordinates:

\bigskip
\begin{theorem}[{\cite[Theorem~2.3]{Ep}}]\label{th:well_known}
A metric tensor \(g : U \to \Sym^2_{\mathrm{reg}} V^*\) defined on a
star-shaped open neighborhood \(U\) of the origin of a vector space
\(V\) is given in geodesic normal coordinates if and only if
\(g(X)_{X,Y} = \langle X, Y\rangle\) for all \(X,Y \in V\). Here,
\(\gg\) denotes the Euclidean structure of \(V\) canonically induced
by \(g\) at the origin.
\end{theorem}

The “only if” direction in Theorem~\ref{th:well_known} is the
classical Gauss lemma. By differentiating the identity
\(g(X)_{X,Y} \equiv \langle X, Y\rangle\) with respect to \(X\) in
\(V\), we see that the \((k+2)\)th coefficient of the Taylor expansion
of the metric tensor in normal coordinates at the origin of \(V\)
belongs to
\begin{equation}\label{eq:formal_Jacobi_operator}
 \scrC_k(V) := \{ h \in \Sym^{k+2} V^* \otimes \Sym^2 V^* \mid
 \forall X,Y \in V : h(X)_{X,Y} = 0 \}
\end{equation}
for all \(k \ge -1\) (where we set \(h(X) := h(X,\ldots,X)\) as
defined earlier). More precisely, we obtain from
Theorem~\ref{th:well_known}:

\bigskip
\begin{corollary}\label{co:well_known}
Let \(V\) be a vector space, \(\gg \in \Sym^2_{\mathrm{reg}} V^*\),
and \(h^j \in \Sym^j V^* \otimes \Sym^2 V^*\) for \(j \ge 1\). The
polynomial
\begin{equation}\label{eq:Taylor2}
 g(X) := \gg + \sum_{j=1}^k h^j(X)
\end{equation}
defines a metric tensor in geodesic normal coordinates on a
sufficiently small star-shaped open neighborhood of the origin if and
only if \(h^j \in \scrC_{j-2}(V)\) for \(j=1,\ldots,k\)
(see~\eqref{eq:formal_Jacobi_operator}).
\end{corollary}

Note that the space \eqref{eq:formal_Jacobi_operator} is nontrivial if
and only if \(k \ge 0\). In fact, the short exact sequence
\begin{equation}\label{eq:ses}
0 \to \scrC_k(V) \longrightarrow \Sym^{k+2} V^* \otimes \Sym^2 V^*
\longrightarrow \Sym^{k+3} V^* \otimes V^* \to 0
\end{equation}
implies that \(\scrC_{-1}(V) = \{0\}\) and that
\[
  \dim \scrC_k(V) = \frac{n(n+1)}{2} \binom{k+n+1}{n-1}
  - n \binom{k+n+2}{n-1}
  = \frac{n(k+1)}{2} \binom{k+n+1}{n-2}
\]
with \(n := \dim(V)\) for \(k \ge 0\), cf. also
Theorem~\ref{th:char2}. Clearly, this formula is in accordance with
Weyl's dimension formula \cite[Theorem~6.3(i)]{FulH} for the dimension
of the irreducible \(\SL_n(\C)\) representation of highest weight
\((k+2,2)\). Furthermore, note that the symmetrized \(k\)th covariant
derivative of the curvature tensor defined in
\eqref{eq:k-ter_Jacobi} satisfies
\(\scrR^k|_p \in \scrC_k(T_pM)\) by the first Bianchi identity.
Therefore, by analogy with Definition~\ref{de:algebraic_curvature_jets},
a collection \(\scrR^{\le k} \in \bigoplus_{j=0}^k \scrC_j(V)\)
should be regarded as an \emph{algebraic} symmetrized \(k\)-jet (of
the curvature tensor associated to some metric \(g\) on \(V\) such
that \(g_0 = \gg\)).

\subsection{Proof of the jet isomorphism theorem in its usual form}
\label{se:fil_aff_spa}

We give a detailed proof of the standard formulation of the jet
isomorphism theorem; for a shorter argument see \cite[p.~77]{FG}.
Let \(M\) be a differentiable manifold and let \(\scrM_{p,k}M\)
denote the system of \(k\)-jets of metric tensors at \(p\)~\cite{W2}.
By definition there are canonical truncation maps
\(\tau_{k,j} : \scrM_{p,k}M \to \scrM_{p,j}M\) for \(j \le k\), and
the isotropy subgroup \(\mathrm{Diff}_p M\) (diffeomorphisms fixing
\(p\)) acts on each \(\scrM_{p,k}M\) on the right.

Two further jet systems play a central role. First, set
\(\bbA_{p,0} M := \bbA_{p,1} M := \Sym^2_{\mathrm{reg}} T^*_p M\)
and, for \(k \ge 0\),
\[
  \bbA_{p,k+2}M := \Sym^2_{\mathrm{reg}} T^*_p M \times
  \bigoplus_{j=0}^k \scrC_j(T_p M)
\]
the space of polynomial metric tensors \(g\) of degree \(\le k+2\)
in geodesic normal coordinates on \(T_p M\), see
Corollary~\ref{co:well_known}. Using an anchored chart
\(f : T_p M \to M\), define the push-forward
\(f_* : \bbA_{p,k}M \to \scrM_{p,k}M\), \(g \mapsto \tilde g := f_* g\).
The induced map
\begin{equation}\label{eq:push_forward}
  f_* : \bbA_{p,k}M \longrightarrow
  \scrM_{p,k}M / \mathrm{Diff}_{p,\Id}M
\end{equation}
is independent of \(f\), where
\(\mathrm{Diff}_{p,\Id}M := \{ f \in \mathrm{Diff}_p M \mid d_p f = \Id \}\).
By the jet isomorphism theorem, this map is an isomorphism.

Second, set \(\bbA_{p,0}^{\star} M := \bbA_{p,1}^{\star} M :=
\Sym^2_{\mathrm{reg}} T^*_p M\) and, for \(k \ge 0\),
\[
  \bbA_{p,k+2}^{\star} M := \Sym^2_{\mathrm{reg}} T^*_p M \times
  \{ \nabla^{\le k}|_0 \rmR \mid \nabla^{\le k}|_0 \rmR
     \text{ is an algebraic \(k\)-jet} \}
\]
the product of regular symmetric bilinear forms with the set of
algebraic curvature \(k\)-jets on \(T_p M\)
(Definition~\ref{de:algebraic_curvature_jets}), together with the
canonical projections that forget higher jet components in the
second factor. By isometry invariance of \(\rmR\) and its covariant
derivatives, evaluation at \(p\) yields a canonical map
\begin{equation}\label{eq:nabla_k}
  \nabla^{\le k} \rmR(M,p) :
  \scrM_{p,k+2}M / \mathrm{Diff}_{p,\Id}M \longrightarrow
  \bbA_{p,k+2}^{\star} M, \quad
  g \longmapsto (g|_p, \nabla^{\le k}|_p \rmR)
\end{equation}
which is again an isomorphism by the theorem. Consequently,
\(\scrM_{p,k+2}M / \mathrm{Diff}_p M \cong
 \bbA_{p,k+2}^{\star} M / \GL(T_p M)\). This is the key point:
any diffeomorphism invariant of the \((k+2)\)-jet of the metric can
be expressed using \(\rmR\) and its covariant derivatives up to order
\(k\), highlighting the centrality of curvature in Riemannian
geometry~\cite{Ep,FG,ABP}. Finally, the proof also shows that the
symmetrization map \(\scrR^{\le k}(T_p M) :
\bbA_{p,k+2}^\star M \to \bbA_{p,k+2}M\), sending
\(\nabla^{\le k}|_p \rmR\) to \(\scrR^{\le k}|_p\), is an isomorphism;
the inverse \(\bbA_{p,k+2}M \to \bbA_{p,k+2}^\star M\) is far from
obvious and identifying it is a main goal of this paper.

The proof uses a canonical affine structure compatible with the
projective structure on \(\bbA_{p,k}M\) and \(\bbA_{p,k}^\star M\).

\bigskip
\begin{definition}[{\cite[p.~60]{KMS}}]\label{de:affine_vector_bundle}
\begin{enumerate}
\item Let \(V\) be a vector space. An affine vector bundle modeled on
\(V\) is a fiber bundle \(\tau : \bbA \to M\) whose fiber \(\tau_p\)
is an affine space modeled on \(V\), and for which the translation
action \(V \times \bbA \to \bbA\) is differentiable. If \(\bbA \to M\)
and \(\bbA' \to M'\) are affine vector bundles modeled on \(V\) and
\(V'\), respectively, a morphism \(F : \bbA \to \bbA'\) is a bundle
morphism such that \(F_{\mathrm{l}}(p_1 - p_2) := F(p_1) - F(p_2)\)
defines a linear map \(F_{\mathrm{l}} : V \to V'\) independent of
\(p_1, p_2 \in \tau_q\) and \(q \in M\), called the associated linear
map.
\item A projective system \((\bbA_{\bullet}, \tau)\) consists of spaces
\(\bbA_{k}\) and maps \(\tau_{j,k} : \bbA_{k} \to \bbA_{j}\) for
\(j \le k\) with \(\tau_{k,k} = \Id_{\bbA_{k}}\) and
\(\tau_{j,k} \circ \tau_{k,\ell} = \tau_{j,\ell}\). Assume there is a
graded vector space \(V_\bullet = \bigoplus_{k = 1}^\infty V_k\) such
that \(\tau_{k-1,k} : \bbA_{k} \to \bbA_{k-1}\) is an affine vector
bundle modeled on \(V_k\) for each \(k \ge 1\). Then
\((\bbA_{\bullet}, \tau)\) is a projective system of affine vector
bundles modeled on \(V_\bullet\). If \((\bbA_{\bullet}, \tau)\) and
\((\bbA_{\bullet}', \tau')\) are modeled on \(V_{\bullet}\) and
\(V_{\bullet}'\), a morphism \(F^{\bullet} : \bbA_{\bullet} \to
\bbA_{\bullet}'\) is a family \(F^i : \bbA_{i} \to \bbA_{i}'\) of
affine maps respecting the projective structure:
\(F^i \circ \tau_{i,j} = \tau'_{i,j} \circ F^j\). The direct sum
\(F^{\bullet}_{\mathrm{agl}} := \bigoplus_{i = 1}^\infty
F^i_{\mathrm{l}} : V_\bullet \to V_\bullet'\) is the associated
graded linear map.
\end{enumerate}
\end{definition}

\bigskip
\begin{remark}\label{re:strategy_for_the_proof}
Given projective systems of affine vector bundles \((\bbA_\bullet, \tau)\)
and \((\bbA_\bullet', \tau')\) modeled on \(V_\bullet\) and
\(V_\bullet'\), an inductive argument shows that
\(F^\bullet : \bbA_\bullet \to \bbA_\bullet'\) is an isomorphism if
and only if
\begin{itemize}
\item \(F^0 : \bbA_0 \to \bbA'_0\) is a diffeomorphism, and
\item the associated graded linear map \(F^\bullet_{\mathrm{agl}} :
V_\bullet \to V_\bullet'\) is an isomorphism.
\end{itemize}
\end{remark}

Let \(V\) be a vector space and write \(\bbA_{\bullet}(V) :=
\bbA_{0,\bullet}V\) and \(\bbA_{\bullet}^\star(V) :=
\bbA_{0,\bullet}^\star V\) for the projective systems of polynomial
metrics in normal coordinates and algebraic curvature jets,
respectively, on \(M := V\) at \(p := 0\). Assuming for the moment
that any algebraic \(k\)-jet \(\nabla^{\le k}|_0 R\) extends to an
algebraic \((k+1)\)-jet \(\nabla^{\le k+1}|_0 R\) (proved in
Theorem~\ref{th:jet_isomorphism_theorem}), the projection maps
\(\bbA_{k}(V) \to \bbA_{j}(V)\) and \(\bbA_{k}^\star(V) \to
\bbA_{j}^\star(V)\) for \(j \le k\) turn \(\bbA_{\bullet}(V)\) and
\(\bbA_{\bullet}^\star(V)\) into projective systems of affine vector
bundles modeled on \(\scrC_{\bullet-2}(V)\) and
\(\scrC_{\bullet-2}^\star(V)\), respectively.

The key maps are the symmetrization \(\scrR^{\le k}(V) :
\bbA^\star_{k+2}(V) \to \bbA_{k+2}(V)\), which is the identity on
\(\Sym^2_{\mathrm{reg}}V^*\) and sends \(\nabla^{\le k}|_0 R\) to
\(\scrR^{\le k}|_0\) (cf.~\eqref{eq:k-ter_Jacobi}), and the
curvature-jet map in the opposite direction \(\nabla^{\le k}\rmR(V) :
\bbA_{k+2}(V) \to \bbA_{k+2}^\star(V)\), which assigns to
\((\gg, h^2, \ldots, h^{k+2}) \in \Sym^2_{\mathrm{reg}}V^* \times
\bigoplus_{j=2}^{k+2} \scrC_{j-2}(V)\) the \(k\)-jet
\(\nabla^{\le k}\rmR|_0\) of the curvature of
\(g(X) := \gg + \sum_{j=2}^{k+2} \frac{1}{j!}\,h^j(X)\) at \(0 \in V\).
The associated graded linear map of \(\scrR^{\le \bullet}(V)\) is
characterized by
\begin{equation}\label{eq:agl_scrS}
  \scrR^k_{\mathrm{agl}}(\nabla^k|_0\rmR) = \scrR^k|_0
\end{equation}
and, as a consequence of the jet isomorphism theorem and Young
symmetrizer theory, the associated graded linear map of
\(\nabla^{\le \bullet}\rmR(V)\) is essentially \(-\tfrac{1}{2}\) times
the Kulkarni–Nomizu product; see
\eqref{eq:agl_nabla_leq_infty_R}. Moreover, by
Theorem~\ref{th:taylor_series} there are noncommutative polynomials
\(Q_k(\scrR^0, \scrR^1, \ldots)\) of degree \(k \ge 2\) giving the
Taylor series of \(g\) in normal coordinates. Set \(Q_0 := \Id\) and
\(Q_1 := 0\), and define
\begin{equation}\label{eq:Q^k}
\begin{array}{l}
\scrQ^k(V) := \bigoplus_{j=0}^k Q_j : \bbA_{k}(V) \longrightarrow
\bbA_{k}(V), \\[0.3em]
\quad (\gg, h^2, \ldots, h^k) \longmapsto (\gg, Q_2(h^2), \ldots,
Q_k(h^2, \ldots, h^k))
\end{array}
\end{equation}

The classical statement of the jet isomorphism theorem from
\cite[Theorem~8.3]{FG}, \cite[Theorem~2.6]{Ep} is as follows.

\bigskip
\begin{theorem}
\label{th:jet_isomorphism_theorem}
\begin{enumerate}
  \item If \((V,\gg)\) is Euclidean, then \(\bbA_{\bullet}^\star(V)\)
  is a projective system of affine vector bundles modeled on
  \(\scrC_{\bullet-2}^\star(V)\). The maps \(\nabla^{\le \bullet} \rmR(V)\),
  \(\scrR^{\le \bullet}(V)\), and \(\scrQ^\bullet(V)\) are
  \(\GL(V)\)-equivariant isomorphisms of such systems, and
  \[
    \scrQ^{\bullet+2}(V) \circ \scrR^{\le \bullet}(V) \circ
    \nabla^{\le \bullet} \rmR(V) = \Id_{\bbA_{\bullet+2}(V)}.
  \]
  \item Let \(M\) be a differentiable manifold, \(p \in M\), and
  \(f : T_pM \to M\) an anchored chart at \(p\). Then
  \(\nabla^{\le \bullet} \rmR(M,p) :
   \scrM_{p,\bullet+2}M / \mathrm{Diff}_{p,\Id}M \to
   \bbA^\star_{p,\bullet+2}M\) (see~\eqref{eq:nabla_k}) and
  \(f_* : \bbA_{p,\bullet}M \to \scrM_{p,\bullet}M /
   \mathrm{Diff}_{p,\Id}M\) (see~\eqref{eq:push_forward})
  are isomorphisms, and
  \[
    f_* \circ \scrQ^{\bullet+2}(T_pM) \circ \scrR^{\le \bullet}(T_pM)
    \circ \nabla^{\le \bullet} \rmR(M,p)
    = \Id_{\scrM_{p,\bullet+2}M / \mathrm{Diff}_{p,\Id}M}.
  \]
\end{enumerate}
\end{theorem}

\begin{proof}
We start with~(a). By Theorem~\ref{th:taylor_series} and
Corollary~\ref{co:well_known}, the composition
\(\scrQ^{k+2}(V) \circ \scrR^{\le k}(V) \circ \nabla^{\le k} \rmR(V)\)
is the identity on \(\bbA_{k+2}(V)\). Moreover,
\(Q_{k+2}(\scrR^0, \scrR^1, \ldots) \equiv
 -2\,\frac{k+1}{k+3}\,\scrR^k\) modulo
\(\{\scrR^0, \ldots, \scrR^{k-1}\}\)
(Section~\ref{se:Taylor_expansion}). Hence the associated graded
linear map \(\scrQ^{\bullet+2}_{\mathrm{agl}}(V)\) is
\begin{equation}\label{eq:agl_scrQ}
  -2 \,\bigoplus_{k=0}^\infty \frac{k+1}{k+3} \,
  \Id_{\scrC_k(V)} \; \longrightarrow\; \scrC_\bullet(V)
\end{equation}
and \(\scrQ^\bullet(V)\) is an isomorphism of projective systems of
affine vector bundles (cf.~Remark~\ref{re:strategy_for_the_proof}).

We prove the following by induction on \(k\):
\begin{itemize}
\item \(\bbA_{k+2}^\star(V)\) is an affine vector bundle over
\(\bbA_{k+1}^\star(V)\) with model \(\scrC^\star_{k}(V)\)
\item \(\scrR^{\le k}(V) : \bbA_{k+2}^\star(V) \to \bbA_{k+2}(V)\) is
an isomorphism of affine vector bundles
\item \(\nabla^{\le k} \rmR(V) : \bbA_{k+2}(V) \to \bbA_{k+2}^\star(V)\)
is an isomorphism of affine vector bundles
\end{itemize}
The claims are obvious for \(k=-2,-1\). For \(k \ge 0\), assume they
hold for \(k-1\). The first is clear for \(k=0\), so let \(k \ge 1\).
Fix \(\gg\) and an algebraic \((k-1)\)-jet \(\nabla^{\le k-1}|_0 \rmR\)
on \(V\). By the induction hypothesis,
\(\nabla^{\le k-1} \rmR(V) : \bbA_{k+1}(V) \to \bbA^\star_{k+1}(V)\)
is an isomorphism, so there exists a metric \(g\) with \(g_0 = \gg\)
whose curvature \((k-1)\)-jet at \(0\) equals \(\nabla^{\le k-1}|_0 \rmR\).
Then \(\nabla^k|_0 \rmR\) extends this to a \(k\)-jet. Hence the fiber
of \(\bbA_{k+2}^\star(V) \to \bbA_{k+1}^\star(V)\) over
\(\nabla^{\le k-1}|_0 \rmR\) is nonempty, and
\(\bbA_{k+2}^\star(V) \to \bbA_{k+1}^\star(V)\) is an affine vector
bundle with model \(\scrC^\star_k(V)\).

Next, the associated linear map of \(\scrR^{\le k}(V)\) is an
isomorphism. Indeed, by Theorems~\ref{th:char1} and \ref{th:char2} in
Appendix~\ref{se:YS} there are alternate realizations
\[
  \scrC_{k}^\star(V) \;=\;
  \bbS^\star_{\scaleto{
  \begin{array}{|c|c|c|c|c|c|}
   \cline{1-5}
    1 & 3 & 5 & \cdots & k+4 \\
   \cline{1-5}
    2 & 4 & \multicolumn{2}{c}{\;\;\;} \\
   \cline{1-2}
  \end{array}}{15pt}} V^*
  \quad\text{and}\quad
  \scrC_k(V) \;=\;
  \bbS_{\scaleto{
  \begin{array}{|c|c|c|c|c|c|}
   \cline{1-5}
    1 & 3 & 5 & \cdots & k+4 \\
   \cline{1-5}
    2 & 4 & \multicolumn{2}{c}{\;\;\;} \\
   \cline{1-2}
  \end{array}}{15pt}} V^*
\]
via Weyl’s construction with the standard Young
tableau~\eqref{eq:Fiedler_alternativ} (Section~\ref{se:YS}). The map
\(\scrC_k^\star(V) \to \scrC_k(V)\) sending
\(\nabla^k|_0 \rmR \mapsto \scrR^k|_0\) equals
\(\frac{-1}{2(k+2)!}\) times the corresponding row symmetrizer, which
is an isomorphism; hence so is \(\scrR^k_{\mathrm{agl}}(V)\) in
\eqref{eq:agl_scrS}. Therefore \(\scrR^{\le k}(V)\) is an isomorphism
of affine vector bundles. Since \(\scrQ^{k+2}(V)\) is an isomorphism
and \(\scrQ^{k+2}(V) \circ \scrR^{\le k}(V) \circ \nabla^{\le k}
\rmR(V) = \Id\) on \(\bbA_{k+2}(V)\), all three maps are
isomorphisms.

In particular, the graded linear map associated with
\(\nabla^{\le \bullet} \rmR(V)\) is
\begin{equation}\label{eq:agl_nabla_leq_infty_R}
  -\frac{1}{2} \owedge \otimes \Id_{\Sym_\bullet V^*} :
  \scrC_{\bullet}(V) \longrightarrow \scrC^\star_{\bullet}(V)
\end{equation}
which completes the proof of~(a).

For~(b), let \(g\) be a metric on \(M\), let \(f : T_pM \to M\) be
geodesic normal coordinates, and set \(\tilde g := f^* g\). By
Theorem~\ref{th:taylor_series},
\(\tilde g = \scrQ^{k+2}(g|_p, \scrR^{\le k}|_p) \mod \scrO(k+3)\).
Hence the \((k+2)\)-jet of \(g = f_* \tilde g\) agrees with
\(f_* \scrQ^{k+2}(g|_p, \scrR^{\le k}|_p)\) up to order \(k+2\), so
\[
  f_* \circ \scrQ^{k+2}(T_pM) \circ \scrR^{\le k}(T_pM) \circ
  \nabla^{\le k} \rmR(M, p)
  = \Id_{\scrM_{p, k+2}M / \mathrm{Diff}_{p, \Id}M}.
\]
On the other hand,
\[
  \nabla^{\le k} \rmR(M, p) \circ f_* =
  \nabla^{\le k} \rmR(T_pM)
\]
by invariance of curvature and its iterated covariant derivatives
under isometries. Therefore
\[
\begin{aligned}
  &\scrQ^{k+2}(T_pM) \circ \scrR^{\le k}(T_pM) \circ
    \nabla^{\le k} \rmR(M, p) \circ f_* \\
  &= \scrQ^{k+2}(T_pM) \circ \scrR^{\le k}(T_pM) \circ
     \nabla^{\le k} \rmR(T_pM) \\
  &= \Id_{\bbA_{k+2}(T_pM)}.
\end{aligned}
\]
Hence both \(f_*\) and \(\nabla^{\le k} \rmR(M, p)\) are
isomorphisms, completing~(b).
\end{proof}

\bigskip
\begin{remark}
It is also possible to obtain a version of the jet isomorphism theorem
for jets of infinite order. For this, one needs a result of
E.~Borel~\cite{BL} which implies that every formal power series is the
Taylor series of some smooth function.
\end{remark}

\section{The inverse of the jet symmetrization map
\texorpdfstring{$\scrR^\bullet$}{R bullet}}
\label{se:Inverse_Jet_Theorem}

Let \((M,g)\) be a Riemannian manifold and \(p\in M\). Our goal is to
reconstruct the \(k\)-jet \(\nabla^{\le k}|_p \rmR\) from the
symmetrized jet \(\scrR^{\le k}|_p\).

\subsection{A direct proof of the Young projection
formula\texorpdfstring{~\eqref{eq:Young_Symmetrizer}}{ (Young
Symmetrizer Eq.)} for linear \texorpdfstring{\(k\)}{k}-jets}
\label{se:Fiedler}

We start with the linear case, i.e.,~we show directly that
\eqref{eq:Young_Symmetrizer} holds for every linear \(k\)-jet
\(\nabla^{\le k}|_p \rmR\).

For \(k=0\) we have to show that
\[
 \rmS^\star_{\scaleto{\begin{array}{|c|c|}
    \cline{1-2}
    1 & 3 \\
    \cline{1-2}
    2 & 4 \\
    \cline{1-2}
   \end{array}
 }{15pt}}\,\rmR_{X_1,X_2,X_3,X_4} \;=\; 12\,\rmR_{X_1,X_2,X_3,X_4}.
\]

\begin{proof}
By definition of the Young symmetrizer,
\begin{align*}
\rmS^\star_{\scaleto{
 \begin{array}{|c|c|}
  \cline{1-2}
  1 & 3 \\
  \cline{1-2}
  2 & 4 \\
  \cline{1-2}
 \end{array}
 }{15pt}}\,\rmR_{X_1,X_2,X_3,X_4}
&= \rmR_{X_1,X_2,X_3,X_4}
 + \rmR_{X_3,X_2,X_1,X_4}
 + \rmR_{X_1,X_4,X_3,X_2}
 + \rmR_{X_3,X_4,X_1,X_2} \\
&\quad - \rmR_{X_2,X_1,X_3,X_4}
 - \rmR_{X_3,X_1,X_2,X_4}
 - \rmR_{X_2,X_4,X_3,X_1}
 - \rmR_{X_3,X_4,X_2,X_1} \\
&\quad - \rmR_{X_1,X_2,X_4,X_3}
 - \rmR_{X_4,X_2,X_1,X_3}
 - \rmR_{X_1,X_3,X_4,X_2}
 - \rmR_{X_4,X_3,X_1,X_2} \\
&\quad + \rmR_{X_2,X_1,X_4,X_3}
 + \rmR_{X_4,X_1,X_2,X_3}
 + \rmR_{X_2,X_3,X_4,X_1}
 + \rmR_{X_4,X_3,X_2,X_1} \\
&= \rmR_{X_1,X_2,X_3,X_4}
 + 2\,\rmR_{X_4,X_3,X_2,X_1}
 + 2\,\rmR_{X_2,X_1,X_4,X_3}
 + 4\,\rmR_{X_3,X_4,X_1,X_2},
\end{align*}
where we used the first Bianchi identity, or equivalently
\(\cyclic_{1,2,3}\rmR_{X_4,X_1,X_2,X_3}=0\). This proves the claim.
\end{proof}

The case \(k=1\) is preceded by the following lemma.

\bigskip
\begin{lemma}\label{le:k=1}
We have
\[
\rmS^\star_{\scaleto{
\begin{array}{|c|c|}
 \cline{1-2}
 1 & 3 \\
 \cline{1-2}
 2 & 4 \\
 \cline{1-2}
\end{array}
}{15pt}}\,\nabla_{X_1}\rmR_{X_3,X_2,X_5,X_4}
\;=\;
\rmS^\star_{\scaleto{
\begin{array}{|c|c|}
 \cline{1-2}
 1 & 3 \\
 \cline{1-2}
 2 & 4 \\
 \cline{1-2}
\end{array}
}{15pt}}\,\nabla_{X_3}\rmR_{X_5,X_2,X_1,X_4}
\;=\; 6\,\nabla_{X_5}\rmR_{X_1,X_2,X_3,X_4},
\]
where the Young symmetrizer acts on the variables
\(X_1,\ldots,X_4\) while \(X_5\) is fixed.
\end{lemma}

\begin{proof}
By pair symmetry,
\[
\rmS^\star_{\scaleto{
\begin{array}{|c|c|}
 \cline{1-2}
 1 & 3 \\
 \cline{1-2}
 2 & 4 \\
 \cline{1-2}
\end{array}
}{15pt}}\,\nabla_{X_1}\rmR_{X_3,X_2,X_5,X_4}
\;=\;
\rmS^\star_{\scaleto{
\begin{array}{|c|c|}
 \cline{1-2}
 1 & 3 \\
 \cline{1-2}
 2 & 4 \\
 \cline{1-2}
\end{array}
}{15pt}}\,\nabla_{X_1}\rmR_{X_5,X_2,X_3,X_4}.
\]
Using the first Bianchi identity,
\begin{align*}
\rmS^\star_{\scaleto{
 \begin{array}{|c|c|}
  \cline{1-2}
  1 & 3 \\
  \cline{1-2}
  2 & 4 \\
  \cline{1-2}
 \end{array}
 }{15pt}}\,\nabla_{X_1}\rmR_{X_3,X_2,X_5,X_4}
&= \nabla_{X_1}\rmR_{X_3,X_2,X_5,X_4}
 + \nabla_{X_3}\rmR_{X_1,X_2,X_5,X_4}
 + \nabla_{X_1}\rmR_{X_3,X_4,X_5,X_2} \\
&\quad + \nabla_{X_3}\rmR_{X_1,X_4,X_5,X_2}
 - \nabla_{X_2}\rmR_{X_3,X_1,X_5,X_4}
 - \nabla_{X_3}\rmR_{X_2,X_1,X_5,X_4} \\
&\quad - \nabla_{X_2}\rmR_{X_3,X_4,X_5,X_1}
 - \nabla_{X_3}\rmR_{X_2,X_4,X_5,X_1}
 - \nabla_{X_1}\rmR_{X_4,X_2,X_5,X_3} \\
&\quad - \nabla_{X_4}\rmR_{X_1,X_2,X_5,X_3}
 - \nabla_{X_1}\rmR_{X_4,X_3,X_5,X_2}
 - \nabla_{X_4}\rmR_{X_1,X_3,X_5,X_2} \\
&\quad + \nabla_{X_2}\rmR_{X_4,X_1,X_5,X_3}
 + \nabla_{X_4}\rmR_{X_2,X_1,X_5,X_3} \\
&\quad + \nabla_{X_2}\rmR_{X_4,X_3,X_5,X_1}
 + \nabla_{X_4}\rmR_{X_2,X_3,X_5,X_1}.
\end{align*}
Using the second Bianchi identity, this equals
\begin{align*}
&\; 3\,\nabla_{X_3}\rmR_{X_1,X_2,X_5,X_4}
 + 3\,\nabla_{X_4}\rmR_{X_2,X_1,X_5,X_3}
 + 3\,\nabla_{X_1}\rmR_{X_3,X_4,X_5,X_2}
 + 3\,\nabla_{X_2}\rmR_{X_4,X_3,X_5,X_1} \\
&= 3\,\nabla_{X_5}\rmR_{X_1,X_2,X_3,X_4}
 + 3\,\nabla_{X_5}\rmR_{X_4,X_3,X_2,X_1}
 = 6\,\nabla_{X_5}\rmR_{X_1,X_2,X_3,X_4}.
\end{align*}
\end{proof}

We are now ready to prove~\eqref{eq:Young_Symmetrizer} for \(k=1\):
\[
 \rmS^\star_{\scaleto{
 \begin{array}{|c|c|c|}
  \cline{1-3}
  1 & 3 & 5\\
  \cline{1-3}
  2 & 4 & \multicolumn{1}{c}{\;\;\;} \\
  \cline{1-2}
 \end{array}}{15pt}}\,\nabla_{X_5}\rmR_{X_1,X_2,X_3,X_4}
 \;=\; 24\,\nabla_{X_5}\rmR_{X_1,X_2,X_3,X_4}.
\]

\begin{proof}
We have
\begin{equation}
 \begin{aligned}
 \rmS^\star_{\scaleto{
 \begin{array}{|c|c|c|}
  \cline{1-3}
  1 & 3 & 5\\
  \cline{1-3}
  2 & 4 & \multicolumn{1}{c}{\;\;\;} \\
  \cline{1-2}
 \end{array}}{15pt}}\,\nabla_{X_5}\rmR_{X_1,X_2,X_3,X_4}
 &= \rmS^\star_{\scaleto{
 \begin{array}{|c|c|}
  \cline{1-2}
  1 & 3 \\
  \cline{1-2}
  2 & 4 \\
  \cline{1-2}
 \end{array}
 }{15pt}}\,\nabla_{X_5}\rmR_{X_1,X_2,X_3,X_4} \\
 &\quad + \rmS^\star_{\scaleto{
 \begin{array}{|c|c|}
  \cline{1-2}
  1 & 3 \\
  \cline{1-2}
  2 & 4 \\
  \cline{1-2}
 \end{array}
 }{15pt}}\,\nabla_{X_1}\rmR_{X_3,X_2,X_5,X_4} \\
 &\quad + \rmS^\star_{\scaleto{
 \begin{array}{|c|c|}
  \cline{1-2}
  1 & 3 \\
  \cline{1-2}
  2 & 4 \\
  \cline{1-2}
 \end{array}
 }{15pt}}\,\nabla_{X_3}\rmR_{X_5,X_2,X_1,X_4}.
 \end{aligned}
\end{equation}
Thus, using the case \(k=0\) together with Lemma~\ref{le:k=1},
\[
 \rmS^\star_{\scaleto{
 \begin{array}{|c|c|c|}
  \cline{1-3}
  1 & 3 & 5\\
  \cline{1-3}
  2 & 4 & \multicolumn{1}{c}{\;\;\;}\\
  \cline{1-2}
 \end{array}
 }{15pt}}\,\nabla_{X_5}\rmR_{X_1,X_2,X_3,X_4}
 \;=\; (12+6+6)\,\nabla_{X_5}\rmR_{X_1,X_2,X_3,X_4}.
\]
\end{proof}

For \(k\ge 2\), we also need the following lemma:

\bigskip
\begin{lemma}\label{le:additional_term1}
For a linear \(2\)-jet \(\nabla^{\le 2}|_p \rmR = (0,0,\nabla^2|_p \rmR)\)
we have
\begin{equation}\label{eq:additional_term1}
 \rmS^\star_{\scaleto{
 \begin{array}{|c|c|}
  \cline{1-2}
  1 & 3 \\
  \cline{1-2}
  2 & 4 \\
  \cline{1-2}
 \end{array}
 }{15pt}}\,\nabla^2_{X_1,X_3}\rmR_{X_5,X_2,X_6,X_4}
 \;=\; 4\,\nabla^{2}_{X_5,X_6}\rmR_{X_1,X_2,X_3,X_4},
\end{equation}
where \(X_5,X_6\) are fixed with respect to the action of the Young
symmetrizer.
\end{lemma}

\begin{proof}
The left-hand side of~\eqref{eq:additional_term1} is
\begin{equation}\label{eq:additional_term11}
 \begin{aligned}
   2\,\nabla^2_{X_1,X_3}\rmR_{X_5,X_2,X_6,X_4}
  &+ 2\,\nabla^2_{X_1,X_3}\rmR_{X_5,X_4,X_6,X_2} \\
  -\,2\,\nabla^2_{X_2,X_3}\rmR_{X_5,X_1,X_6,X_4}
  &- 2\,\nabla^2_{X_2,X_3}\rmR_{X_5,X_4,X_6,X_1} \\
  -\,2\,\nabla^2_{X_1,X_4}\rmR_{X_5,X_2,X_6,X_3}
  &- 2\,\nabla^2_{X_1,X_4}\rmR_{X_5,X_3,X_6,X_2} \\
  +\,2\,\nabla^2_{X_2,X_4}\rmR_{X_5,X_1,X_6,X_3}
  &+ 2\,\nabla^2_{X_2,X_4}\rmR_{X_5,X_3,X_6,X_1},
 \end{aligned}
\end{equation}
where we used the trivial Ricci identity
\(\nabla^2_{X,Y}\rmR=\nabla^2_{Y,X}\rmR\).

Applying the second Bianchi identity together with pair symmetry to
each pair of summands occupying the same position in lines one and
three, or two and four, and using other curvature symmetries, this
becomes
\[
 2\,\rmS_{\scaleto{
  \begin{array}{|c|c|}
   \cline{1-2}
   5 & 6 \\
   \cline{1-2}
  \end{array}}{7.5pt}}\,
 \rmS_{\scaleto{
  \begin{array}{|c|}
   \cline{1-1}
   1 \\
   \cline{1-1}
   2 \\
   \cline{1-1}
  \end{array}}{15pt}}\,\nabla^2_{X_1,X_5}\rmR_{X_6,X_2,X_3,X_4}.
\]
Using again the trivial Ricci identity and the second Bianchi identity,
\begin{align*}
 \rmS_{\scaleto{
  \begin{array}{|c|}
   \cline{1-1}
   1 \\
   \cline{1-1}
   2 \\
   \cline{1-1}
  \end{array}}{15pt}}\,\nabla^2_{X_1,X_5}\rmR_{X_6,X_2,X_3,X_4}
 &= \rmS_{\scaleto{
  \begin{array}{|c|}
   \cline{1-1}
   1 \\
   \cline{1-1}
   2 \\
   \cline{1-1}
  \end{array}}{15pt}}\,\nabla^2_{X_5,X_1}\rmR_{X_6,X_2,X_3,X_4} \\
 &= \nabla^2_{X_5,X_6}\rmR_{X_1,X_2,X_3,X_4}.
\end{align*}
Using the trivial Ricci identity once more yields the claimed result.
\end{proof}

\paragraph{Proof of~\eqref{eq:Young_Symmetrizer} for \(k\ge 2\)}
Suppose that \(\nabla^{\le k}|_p \rmR\) is a linear \(k\)-jet,
i.e.,~\(\nabla^{\le k}|_p \rmR=(0,\ldots,0,\nabla^k|_p \rmR)\).
The natural right action of the symmetric group
\(S_{\{1,3,5,\ldots,k+4\}}\) on
\(\nabla^k_{X_5,\ldots,X_{k+4}}\rmR_{X_1,X_2,X_3,X_4}\)
factorizes over the space of right cosets
\[
S_{\{1,3,5,\ldots,k+4\}} \,/\, S_{\{5,\ldots,k+4\}}.
\]
To find a suitable set of representatives, note that there is a
canonical inclusion
\(S_{\{1,3\}} \hookrightarrow
 S_{\{1,3,5,\ldots,k+4\}} \,/\, S_{\{5,\ldots,k+4\}}\)
yielding two right cosets. Similarly, for each \(A=5,\ldots,k+4\)
there is a natural inclusion
\(S_{\{1,3,A\}} \hookrightarrow
 S_{\{1,3,5,\ldots,k+4\}} \,/\, S_{\{5,\ldots,k+4\}}\),
which produces
\(|S_{\{1,3,A\}} \setminus S_{\{1,3\}}| = 6-2 = 4\) further right
cosets. Together with the \((k-1)k\) permutations \((1\,A)(3\,B)\) for
\(5\le A\ne B\le k+4\), we obtain
\[
2 + 4k + (k-1)k = k^2 + 3k + 2 = (k+1)(k+2)
\]
distinct elements of \(S_{\{1,3,5,\ldots,k+4\}}\) that exhaust the
space of right cosets.

In the free vector space over
\(S_{\{1,3,5,\ldots,k+4\}} \,/\, S_{\{5,\ldots,k+4\}}\),
\begin{equation*}
 \sum_{[\pi] \in S_{\{1,3,5,\ldots,k+4\}} / S_{\{5,\ldots,k+4\}}}
 [\pi]
 = \sum_{A=5}^{k+4} \sum_{\pi \in S_{\{1,3,A\}}} [\pi]
   - (k-1) \sum_{\pi \in S_{\{1,3\}}} [\pi]
   + \sum_{A \neq B = 5}^{k+4} [(1\,A)(3\,B)].
\end{equation*}

Hence, using the cases \(k=0\) and \(k=1\) together with
Lemma~\ref{le:additional_term1}, we have
\begin{align*}
 \frac{1}{k!}\,
 \rmS^\star_{\scaleto{
 \begin{array}{|c|c|c|c|c|}
  \cline{1-4}
  1 & 3 & \cdots & k+4\\
  \cline{1-4}
  2 & 4 & \multicolumn{2}{c}{\;} \\
  \cline{1-2}
 \end{array}
 }{15pt}}\,
 \nabla^k_{X_5,\ldots,X_{k+4}}\rmR_{X_1,X_2,X_3,X_4}
&= \bigg(
   \sum_{A=5}^{k+4}
   \rmS^\star_{\scaleto{
   \begin{array}{|c|c|c|c|}
    \cline{1-3}
    1 & 3 & A\\
    \cline{1-3}
    2 & 4 & \multicolumn{1}{c}{\;} \\
    \cline{1-2}
   \end{array}
   }{15pt}}
   -(k-1)\,
   \rmS^\star_{\scaleto{
   \begin{array}{|c|c|}
    \cline{1-2}
    1 & 3 \\
    \cline{1-2}
    2 & 4 \\
    \cline{1-2}
   \end{array}
   }{15pt}}
   \bigg)
   \nabla^k_{X_5,\ldots,X_{k+4}}\rmR_{X_1,X_2,X_3,X_4} \\
&\quad
   + \sum_{5=A<B}^{k+4}
   \rmS^\star_{\scaleto{
   \begin{array}{|c|c|}
    \cline{1-2}
    1 & 3 \\
    \cline{1-2}
    2 & 4 \\
    \cline{1-2}
   \end{array}
   }{15pt}}
   \nabla^k_{X_1,X_3,X_5,\ldots,\hat{X}_A,\ldots,\hat{X}_B,
   \ldots,X_{k+4}}
   \rmR_{X_A,X_2,X_B,X_4} \\
&= \underbrace{ \big(24k - 12(k-1) + 4\,\tfrac{k(k-1)}{2}\big)
   }_{=\,2(k+2)(k+3)}
   \nabla^k_{X_5,\ldots,X_{k+4}}\rmR_{X_1,X_2,X_3,X_4}.
\end{align*}
Since \(h_k = 2\,k!\,(k+2)(k+3)\), this proves
\eqref{eq:Young_Symmetrizer} for all \(k\ge 0\). \qed

\subsection{Generalization of\texorpdfstring{~\eqref{eq:Young_Symmetrizer}}{ (Young Symmetrizer Eq.)} for arbitrary \texorpdfstring{\(k\)}{k}-jets}
\label{se:general_case}

We now reconstruct the \(k\)-jet \(\nabla^{\le k}|_p \rmR\) from its
symmetrization \(\scrR^{\le k}|_p\) for an arbitrary Riemannian
manifold. In this general case, the Ricci identity in the
form~\eqref{eq:Ricci_R} must also be taken into account.
Lemma~\ref{le:additional_term1} then admits the following
modification.

\bigskip
\begin{lemma}\label{le:additional_term2}
For an arbitrary Riemannian manifold, we have
\begin{equation}\label{eq:additional_term2}
\begin{array}{lll}
\rmS^\star_{\scaleto{\begin{array}{|c|c|}
\cline{1-2}
1 & 3 \\
\cline{1-2}
2 & 4 \\
\cline{1-2}
\end{array}}{15pt}}\,
\nabla^2_{X_1,X_3}\rmR_{X_2,X_5,X_4,X_6}
&=&
\rmS_{\scaleto{\begin{array}{|c|c|}
\cline{1-2}
5 & 6 \\
\cline{1-2}
\end{array}}{7.5pt}}\!
\big(
  2\,\nabla^2_{X_5,X_6}\rmR_{X_1,X_2,X_3,X_4}
\\&&\quad
  +\,2\,
  \rmS_{\scaleto{\begin{array}{|c|}
  \cline{1-1}
  1 \\
  \cline{1-1}
  2 \\
  \cline{1-1}
  \end{array}}{15pt}}\,
  \rmR_{X_1,X_5}\,\rmR_{X_6,X_2,X_3,X_4}
\\&&\quad
  -\,
  \rmS_{\scaleto{\begin{array}{|c|}
  \cline{1-1}
  1 \\
  \cline{1-1}
  2 \\
  \cline{1-1}
  \end{array}}{15pt}}\,
  \rmS_{\scaleto{\begin{array}{|c|}
  \cline{1-1}
  3 \\
  \cline{1-1}
  4 \\
  \cline{1-1}
  \end{array}}{15pt}}\,
  \rmR_{X_1,X_3}\,\rmR_{X_5,X_2,X_6,X_4}
\big)
\end{array}
\end{equation}
\end{lemma}

\begin{proof}
The argument parallels the proof of Lemma~\ref{le:additional_term1},
but with the Ricci identity contributing additional curvature terms.
Starting from the analogue of~\eqref{eq:additional_term11} and adding
the curvature term obtained from the Ricci identity gives
\begin{align*}
&\rmS^\star_{\scaleto{\begin{array}{|c|c|}
\hline 1 & 3 \\ \hline 2 & 4 \\ \hline
\end{array}}{15pt}}
\,\nabla^2_{X_1,X_3}\rmR_{X_2,X_5,X_4,X_6} \\
&\quad=
\rmS_{\scaleto{\begin{array}{|c|c|}
\hline 5 & 6 \\ \hline
\end{array}}{7.5pt}}
\big(
  2\,
  \rmS_{\scaleto{\begin{array}{|c|}
  \hline 1 \\ \hline 2 \\ \hline
  \end{array}}{15pt}}
  \nabla^2_{X_1,X_5}\rmR_{X_6,X_2,X_3,X_4}
  -
  \rmS_{\scaleto{\begin{array}{|c|}
  \hline 1 \\ \hline 2 \\ \hline
  \end{array}}{15pt}}
  \rmS_{\scaleto{\begin{array}{|c|}
  \hline 3 \\ \hline 4 \\ \hline
  \end{array}}{15pt}}
  \rmR_{X_1,X_3}\,\rmR_{X_5,X_2,X_6,X_4}
\big)
\end{align*}
Continuing exactly as in Lemma~\ref{le:additional_term1} and applying
the Ricci identity a second time yields~\eqref{eq:additional_term2}.
\end{proof}

To find the right modification of~\eqref{eq:Young_Symmetrizer}, let
\((\bbE,\nabla^\bbE)\) be a vector bundle with a linear connection
over \((M,g)\) (e.g., a tensor bundle with the induced connection).
Following~\cite[Ch.~4]{W1}, the symmetrized iterated \(k\)th covariant
derivative of a section \(\psi \in \Gamma(\bbE)\) is defined by
\begin{equation}\label{eq:k_jet}
\jet^k_{X_1,\ldots,X_k}\psi
:= \frac{1}{k!} \sum_{\sigma \in \rmS_k}
\nabla^k_{X_{\sigma(1)},\ldots,X_{\sigma(k)}}\psi .
\end{equation}
We obtain the following modification of the Young projection
formula~\eqref{eq:Young_Symmetrizer}.

\bigskip
\begin{proposition}\label{p:Young_Symmetrizer_modifiziert_1}
Let \((\rmR|_p,\nabla|_p \rmR,\ldots,\nabla^k|_p \rmR)\) be an
arbitrary curvature \(k\)-jet. The term
\begin{equation*}
\Big(
 \frac{1}{k!}\,
 \rmS^\star_{\scaleto{
 \begin{array}{|c|c|c|c|c|c|}
 \cline{1-5}
 1 & 3 & 5 & \cdots & k+4 \\
 \cline{1-5}
 2 & 4 & \multicolumn{1}{c}{\;} \\
 \cline{1-2}
 \end{array}}{15pt}}
 \nabla^k_{X_5,\ldots,X_{k+4}}
 - 2(k+2)(k+3)\,\jet^k_{X_5,\ldots,X_{k+4}}
\Big)
\rmR_{X_1,X_2,X_3,X_4}
\end{equation*}
is given by~\eqref{eq:Young_Symmetrizer_modifiziert_1}.
\end{proposition}

\begin{proof}
As in the proof of~\eqref{eq:Young_Symmetrizer} at the end of
Section~\ref{se:Fiedler}, we have
\begin{align*}
&\frac{1}{k!}\,
\rmS^\star_{\scaleto{
\begin{array}{|c|c|c|c|c|}
\cline{1-4}
1 & 3 & \cdots & k+4 \\
\cline{1-4}
2 & 4 & \multicolumn{2}{c}{\;} \\
\cline{1-2}
\end{array}}{15pt}}
\nabla^k_{X_5,\ldots,X_{k+4}}\,
\rmR_{X_1,X_2,X_3,X_4} \\
&\quad=
\Big(
 \sum_{A=1}^k
 \rmS^\star_{\scaleto{
 \begin{array}{|c|c|c|}
 \cline{1-3}
 1 & 3 & A+4 \\
 \cline{1-3}
 2 & 4 & \multicolumn{1}{c}{\;} \\
 \cline{1-2}
 \end{array}}{15pt}}
 -(k-1)\,
 \rmS^\star_{\scaleto{
 \begin{array}{|c|c|}
 \cline{1-2}
 1 & 3 \\
 \cline{1-2}
 2 & 4 \\
 \cline{1-2}
 \end{array}}{15pt}}
\Big)
\jet^k_{X_5,\ldots,X_{k+4}}\,
\rmR_{X_1,X_2,X_3,X_4} \\
&\qquad+
\sum_{\substack{A,B=1\\A<B}}^k
\rmS^\star_{\scaleto{
\begin{array}{|c|c|}
\cline{1-2}
1 & 3 \\
\cline{1-2}
2 & 4 \\
\cline{1-2}
\end{array}}{15pt}}
\jet^k_{X_1,X_3,X_5,\ldots,\hat{X}_{A+4},\ldots,
\hat{X}_{B+4},\ldots,X_{k+4}}\,
\rmR_{X_{A+4},X_2,X_{B+4},X_4} .
\end{align*}

Furthermore, writing
\begin{align*}
&\rmS^\star_{\scaleto{
\begin{array}{|c|c|c|c|}
\cline{1-3}
1 & 3 & A+4 \\
\cline{1-3}
2 & 4 & \multicolumn{1}{c}{\;} \\
\cline{1-2}
\end{array}}{15pt}}
\jet^k_{X_5,\ldots,X_{k+4}}\,
\rmR_{X_1,X_2,X_3,X_4} \\
&\quad=
\rmS^\star_{\scaleto{
\begin{array}{|c|c|c|c|}
\cline{1-3}
1 & 3 & A+4 \\
\cline{1-3}
2 & 4 & \multicolumn{1}{c}{\;} \\
\cline{1-2}
\end{array}}{15pt}}
\jet^{k-1}_{X_5,\ldots,\hat{X}_{A+4},\ldots,X_{k+4}}\,
\nabla_{X_{A+4}}\rmR_{X_1,X_2,X_3,X_4} \\
&\qquad+
\rmS^\star_{\scaleto{
\begin{array}{|c|c|c|c|}
\cline{1-3}
1 & 3 & A+4 \\
\cline{1-3}
2 & 4 & \multicolumn{1}{c}{\;} \\
\cline{1-2}
\end{array}}{15pt}}
\Big(
 \jet^k_{X_5,\ldots,X_{k+4}}
 -\jet^{k-1}_{X_5,\ldots,\hat{X}_{A+4},\ldots,X_{k+4}}\,
  \nabla_{X_{A+4}}
\Big)
\rmR_{X_1,X_2,X_3,X_4} ,
\end{align*}
we see that
\begin{align*}
&\sum_{A=1}^k
\rmS^\star_{\scaleto{
\begin{array}{|c|c|c|c|}
\cline{1-3}
1 & 3 & A+4 \\
\cline{1-3}
2 & 4 & \multicolumn{1}{c}{\;} \\
\cline{1-2}
\end{array}}{15pt}}
\jet^k_{X_5,\ldots,X_{k+4}}\,
\rmR_{X_1,X_2,X_3,X_4} \\
&\quad=
24\,k\,\jet^k_{X_5,\ldots,X_{k+4}}\,
\rmR_{X_1,X_2,X_3,X_4} \\
&\qquad+
\sum_{A=1}^k
\rmS^\star_{\scaleto{
\begin{array}{|c|c|c|c|}
\cline{1-3}
1 & 3 & A+4 \\
\cline{1-3}
2 & 4 & \multicolumn{1}{c}{\;} \\
\cline{1-2}
\end{array}}{15pt}}
\Big(
 \jet^k_{X_5,\ldots,X_{k+4}}
 -\jet^{k-1}_{X_5,\ldots,\hat{X}_{A+4},\ldots,X_{k+4}}\,
  \nabla_{X_{A+4}}
\Big)
\rmR_{X_1,X_2,X_3,X_4} ,
\end{align*}
because \(h_1=24\).

Next,
\[
\rmS^\star_{\scaleto{
\begin{array}{|c|c|}
\cline{1-2}
1 & 3 \\
\cline{1-2}
2 & 4 \\
\cline{1-2}
\end{array}}{15pt}}
\jet^k_{X_5,\ldots,X_{k+4}}\,
\rmR_{X_1,X_2,X_3,X_4}
= 12\,\jet^k_{X_5,\ldots,X_{k+4}}\,
\rmR_{X_1,X_2,X_3,X_4}
\]
since \(h_0=12\).

Also, using Lemma~\ref{le:additional_term2},
\begin{align*}
&\sum_{\substack{A,B=1\\A<B}}^k
\rmS^\star_{\scaleto{
\begin{array}{|c|c|}
\cline{1-2}
1 & 3 \\
\cline{1-2}
2 & 4 \\
\cline{1-2}
\end{array}}{15pt}}
\jet^k_{X_1,X_3,X_5,\ldots,\hat{X}_{A+4},\ldots,
\hat{X}_{B+4},\ldots,X_{k+4}}\,
\rmR_{X_{A+4},X_2,X_{B+4},X_4} \\
&\quad=
2k(k-1)\,\jet^k_{X_5,\ldots,X_{k+4}}\,
\rmR_{X_1,X_2,X_3,X_4} \\
&\qquad+
\sum_{\substack{A,B=1\\A\neq B}}^k
\jet^{k-2}_{X_5,\ldots,\hat{X}_{A+4},\ldots,
\hat{X}_{B+4},\ldots,X_{k+4}}
\Big(
 2\,
 \rmS_{\scaleto{
 \begin{array}{|c|}
 \cline{1-1}
 1 \\
 \cline{1-1}
 2 \\
 \cline{1-1}
 \end{array}}{15pt}}\,
 \rmR_{X_1,X_{A+4}}\,
 \rmR_{X_{B+4},X_2,X_3,X_4} \\
&\qquad\qquad
 -\,
 \rmS_{\scaleto{
 \begin{array}{|c|}
 \cline{1-1}
 1 \\
 \cline{1-1}
 2 \\
 \cline{1-1}
 \end{array}}{15pt}}
 \rmS_{\scaleto{
 \begin{array}{|c|}
 \cline{1-1}
 3 \\
 \cline{1-1}
 4 \\
 \cline{1-1}
 \end{array}}{15pt}}\,
 \rmR_{X_1,X_3}\,
 \rmR_{X_{A+4},X_2,X_{B+4},X_4}
\Big) \\
&\qquad+
\sum_{\substack{A,B=1\\A<B}}^k
\rmS^\star_{\scaleto{
\begin{array}{|c|c|}
\cline{1-2}
1 & 3 \\
\cline{1-2}
2 & 4 \\
\cline{1-2}
\end{array}}{15pt}}
\Big(
 \jet^k_{X_1,X_3,X_5,\ldots,\hat{X}_{A+4},\ldots,
 \hat{X}_{B+4},\ldots,X_{k+4}}
 -\,
 \jet^{k-2}_{X_5,\ldots,\hat{X}_{A+4},\ldots,
 \hat{X}_{B+4},\ldots,X_{k+4}}\,
 \nabla^2_{X_1,X_3}
\Big)\!
\rmR_{X_{A+4},X_2,X_{B+4},X_4} .
\end{align*}
Using these considerations, we obtain the desired result, in analogy
with the proof of~\eqref{eq:Young_Symmetrizer}.
\end{proof}

\subsection{Formulas relating \texorpdfstring{$\jet^k\psi$}{jet\^k psi}
and \texorpdfstring{$\nabla^k\psi$}{nabla\^k psi}}
\label{se:jet_theory}

Let \((\bbE,\nabla)\) be a vector bundle with a linear connection
over \((M,g)\), and let \(\psi\) be a section of \(\bbE\). The
goal of this section is to relate the symmetrized iterated covariant
derivative \(\jet^k\psi\) to the iterated covariant derivative
\(\nabla^k\psi\) itself; see Proposition~\ref{p:Simple_Jet-Formula}. In
the following, we use the Ricci identity in the form stated
in~\eqref{eq:Ricci_psi}, i.e.,~the Leibniz rule is not yet
incorporated. We then have the following simple jet formula.

\bigskip
\begin{lemma}\label{le:Simple_Jet-Formula}
Let \(\bbE\) be a vector bundle over \((M,g)\) equipped with a linear
connection \(\nabla\). For every \(\psi \in \Gamma(\bbE)\),
\begin{equation}\label{eq:Simple_Jet-Formula}
 \nabla^k_{X, \ldots, X, Y}\,\psi - \jet^{k}_{X, \ldots, X, Y}\,\psi
 = \frac{1}{k} \sum_{j=1}^{k-1} j\,\nabla^{j-1}_{X, \ldots, X}\,\rmR_{X, Y}\,
 \nabla^{k-j-1}_{X, \ldots, X}\psi .
\end{equation}
\end{lemma}

\begin{proof}
Using a telescopic sum argument and the Ricci
identity~\eqref{eq:Ricci_psi},
\begin{align*}
 &\big(\nabla^k_{X, \ldots, X, \underset{i}{Y}, X, \ldots, X}
 - \nabla^k_{X, \ldots, X, Y}\big)\,\psi \\
 &\quad= \sum_{j=i}^{k-1}
 \big(\nabla^k_{X, \ldots, X, \underset{j}{Y}, X, \ldots, X}
 - \nabla^k_{X, \ldots, X, \underset{j+1}{Y}, X, \ldots, X}\big)\,\psi \\
 &\quad= \sum_{j=i}^{k-1} \nabla^{j-1}_{X, \ldots, X}\,\rmR_{Y, X}\,
 \nabla^{k-j-1}_{X, \ldots, X}\psi .
\end{align*}
From this, it follows that
\[
 k\,\nabla^k_{X, \ldots, X, Y}\,\psi
 - \sum_{i=1}^k \nabla_{X, \ldots, X, \underset{i}{Y}, X, \ldots, X}\psi
 = \sum_{j=1}^{k-1} j\,\nabla^{j-1}_{X, \ldots, X}\,\rmR_{X, Y}\,
 \nabla^{k-j-1}_{X, \ldots, X}\psi ,
\]
which gives~\eqref{eq:Simple_Jet-Formula} after dividing by \(k\)
(using \(\rmR_{Y,X}=-\rmR_{X,Y}\)).
\end{proof}

We now generalize~\eqref{eq:Simple_Jet-Formula} to obtain a formula
for the difference \(\jet^{k}\psi - \jet^{\ell}\nabla^{k-\ell}\psi\)
for \(0 \le \ell \le k\). For this, we define for every
\(\sigma \in \rmS_k\) different from \(\Id_{\{1,\ldots,k\}}\) the
number \(k_\sigma\) to be the largest index such that
\(\sigma(k_\sigma) \neq k_\sigma\). In other words, with respect to the
canonical inclusion \(\rmS_\ell \subset \rmS_k\), we have
\(\sigma \in \rmS_k \setminus \rmS_\ell\) iff \(k_\sigma \ge \ell+1\)
(where \(\setminus\) denotes the relative complement). Moreover, put
\(j_\sigma := \sigma^{-1}(k_\sigma)\). In the same notation as in
Lemma~\ref{le:Simple_Jet-Formula} we have:

\bigskip
\begin{proposition}\label{p:Simple_Jet-Formula}
Let \(\bbE\) be a vector bundle over \((M,g)\) equipped with a linear
connection \(\nabla\). For every \(\psi \in \Gamma(\bbE)\) and all
\(0 \le \ell \le k\),
\begin{equation}\label{eq:jet^k-jet^l}
\begin{aligned}
 &\big(\jet^{k}_{X_1, \ldots, X_k} - \jet^{\ell}_{X_1, \ldots, X_\ell}
 \nabla^{k-\ell}_{X_{\ell+1}, \ldots, X_k}\big)\psi \\
 &\quad=
 \sum_{\sigma \in \rmS_k \setminus \rmS_\ell} \frac{j_\sigma}{k_\sigma!}\,
 \nabla^{j_\sigma-1}_{X_{\sigma(1)}, \ldots, X_{\sigma(j_\sigma-1)}}\,
 \rmR_{X_{k_\sigma}, X_{\sigma(j_\sigma+1)}}\,
 \nabla^{k-j_\sigma-1}_{X_{\sigma(j_\sigma+2)}, \ldots, X_{\sigma(k)}}\psi .
\end{aligned}
\end{equation}
In particular,
\begin{equation}\label{eq:jet^k-nabla^k}
\begin{aligned}
 &\jet^{k}_{X_1, \ldots, X_k}\psi - \nabla^k_{X_1, \ldots, X_k}\psi \\
 &\quad=
 \sum_{\substack{\sigma \in \rmS_k \\ \sigma \ne \Id}}\frac{j_\sigma}{k_\sigma!}\,
 \nabla^{j_\sigma-1}_{X_{\sigma(1)}, \ldots, X_{\sigma(j_\sigma-1)}}\,
 \rmR_{X_{k_\sigma}, X_{\sigma(j_\sigma+1)}}\,
 \nabla^{k-j_\sigma-1}_{X_{\sigma(j_\sigma+2)}, \ldots, X_{\sigma(k)}}\psi .
\end{aligned}
\end{equation}
\end{proposition}

\begin{proof}
We proceed by induction on \(k\). For \(k=0\) there is nothing to
show. Assume the claim holds for all vector bundles and some integer
\(k \ge 0\). We prove it also holds for \(k+1\). For \(\ell = k+1\)
there is again nothing to show, so we may assume \(\ell \le k\).

Applying the induction hypothesis to the vector bundle
\(\bbE \otimes T^*M\) with the induced connection (again denoted by
\(\nabla\)) and the section \(\nabla\psi\) of this vector bundle, we
obtain
\begin{align*}
&\big(\jet^{k}_{X_1, \ldots, X_k}\,\nabla_{X_{k+1}}
 - \jet^{\ell}_{X_1, \ldots, X_\ell}
   \nabla^{k+1-\ell}_{X_{\ell+1}, \ldots, X_{k+1}}\big)\psi \\
&\quad=
\sum_{\sigma \in \rmS_k \setminus \rmS_\ell} \frac{j_\sigma}{k_\sigma!}\,
\nabla^{j_\sigma-1}_{X_{\sigma(1)}, \ldots, X_{\sigma(j_\sigma-1)}}\,
\rmR_{X_{k_\sigma}, X_{\sigma(j_\sigma+1)}}\,
\nabla^{k-j_\sigma-1}_{X_{\sigma(j_\sigma+2)}, \ldots, X_{\sigma(k)}}\,
\nabla_{X_{k+1}}\psi .
\end{align*}

Furthermore, polarizing~\eqref{eq:Simple_Jet-Formula} yields
\begin{align*}
&\big(\jet^{k+1}_{X_1, \ldots, X_k, X_{k+1}}
 - \jet^k_{X_1, \ldots, X_k}\,\nabla_{X_{k+1}}\big)\psi \\
&\quad=
\sum_{j=1}^k
\sum_{\substack{\sigma \in \rmS_{k+1} \\ \sigma(j)=k+1}}
\frac{j}{(k+1)!}\,
\nabla^{j-1}_{X_{\sigma(1)}, \ldots, X_{\sigma(j-1)}}\,
\rmR_{X_{k+1}, X_{\sigma(j+1)}}\,
\nabla^{k-j}_{X_{\sigma(j+2)}, \ldots, X_{\sigma(k+1)}}\psi \\
&\quad=
\sum_{\sigma \in \rmS_{k+1} \setminus \rmS_k} \frac{j_\sigma}{k_\sigma!}\,
\nabla^{j_\sigma-1}_{X_{\sigma(1)}, \ldots, X_{\sigma(j_\sigma-1)}}\,
\rmR_{X_{k+1}, X_{\sigma(j_\sigma+1)}}\,
\nabla^{k-j_\sigma}_{X_{\sigma(j_\sigma+2)}, \ldots, X_{\sigma(k+1)}}\psi ,
\end{align*}
where we used for the second equality that \(\sigma^{-1}(k+1) \le k\)
holds iff \(\sigma \in \rmS_{k+1} \setminus \rmS_k\) and that
\(k_\sigma = k+1\) for such \(\sigma\).

Using a telescopic sum, we therefore have
\begin{align*}
&\big(\jet^{k+1}_{X_1, \ldots, X_{k+1}}
 - \jet^{\ell}_{X_1, \ldots, X_\ell}
   \nabla^{k+1-\ell}_{X_{\ell+1}, \ldots, X_{k+1}}\big)\psi \\
&\quad=
\Big(
\big(\jet^{k+1}_{X_1, \ldots, X_{k+1}}
 - \jet^k_{X_1, \ldots, X_k}\,\nabla_{X_{k+1}}\big) \\
&\qquad
+ \big(\jet^k_{X_1, \ldots, X_k}\,\nabla_{X_{k+1}}
 - \jet^{\ell}_{X_1, \ldots, X_\ell}
   \nabla^{k+1-\ell}_{X_{\ell+1}, \ldots, X_{k+1}}\big)
\Big)\psi \\
&\quad=
\sum_{\sigma \in \rmS_{k+1} \setminus \rmS_k} \frac{j_\sigma}{k_\sigma!}\,
\nabla^{j_\sigma-1}_{X_{\sigma(1)}, \ldots, X_{\sigma(j_\sigma-1)}}\,
\rmR_{X_{k+1}, X_{\sigma(j_\sigma+1)}}\,
\nabla^{k-j_\sigma}_{X_{\sigma(j_\sigma+2)}, \ldots, X_{\sigma(k+1)}}\psi \\
&\qquad
+ \sum_{\sigma \in \rmS_k \setminus \rmS_\ell} \frac{j_\sigma}{k_\sigma!}\,
\nabla^{j_\sigma-1}_{X_{\sigma(1)}, \ldots, X_{\sigma(j_\sigma-1)}}\,
\rmR_{X_{k_\sigma}, X_{\sigma(j_\sigma+1)}}\,
\nabla^{k-j_\sigma-1}_{X_{\sigma(j_\sigma+2)}, \ldots, X_{\sigma(k)}}\,
\nabla_{X_{k+1}}\psi \\
&\quad=
\sum_{\sigma \in \rmS_{k+1} \setminus \rmS_\ell} \frac{j_\sigma}{k_\sigma!}\,
\nabla^{j_\sigma-1}_{X_{\sigma(1)}, \ldots, X_{\sigma(j_\sigma-1)}}\,
\rmR_{X_{k_\sigma}, X_{\sigma(j_\sigma+1)}}\,
\nabla^{k-j_\sigma}_{X_{\sigma(j_\sigma+2)}, \ldots, X_{\sigma(k+1)}}\psi ,
\end{align*}
where we used the decomposition
\(\rmS_{k+1} \setminus \rmS_\ell =
 (\rmS_{k+1} \setminus \rmS_k) \dot{\cup} (\rmS_k \setminus \rmS_\ell)\).

This completes the induction step for \(k+1\). By setting \(\ell := 0\)
we obtain the claimed formula for \(\nabla^k_{X_1, \ldots, X_k}\psi\).
\end{proof}

For a version of Lemma~\ref{le:Simple_Jet-Formula} and
Proposition~\ref{p:Simple_Jet-Formula} with the Leibniz rule
incorporated, one simply uses \eqref{eq:Ricci_refined} instead of
\eqref{eq:Ricci_psi}.

\bigskip
\begin{corollary}\label{co:B}
Let \((M,g)\) be an arbitrary Riemannian manifold with Levi–Civita
connection \(\nabla\) and curvature tensor \(\rmR\). For each \(k \ge 0\)
there exists a quadratic expression \(B(\nabla^{\le k-2}\rmR)\) in the
\((k-2)\)-jet such that~\eqref{eq:YS2} holds.
\end{corollary}

\begin{proof}
By Proposition~\ref{p:Simple_Jet-Formula} (applied to the vector
bundle \(\bbE := \scrC_0^\star(TM)\) of algebraic curvature tensors,
the section \(\psi := \rmR\), and the connection
\(\nabla^{\scrC_0^\star TM}\) induced by the Levi–Civita connection),
together with the Leibniz rule expressed in~\eqref{eq:Ricci_refined},
it follows that~\eqref{eq:Young_Symmetrizer_modifiziert_1} is a
quadratic expression in \(\nabla^{\le k-2}\rmR\).

Also \((\nabla^k - \jet^k)\rmR\) is a quadratic term
\(B_1(\nabla^{\le k-2}\rmR)\) in the \((k-2)\)-jet. Then we have
\begin{align}
 \nabla^k \rmR + \frac{k+1}{k+3}(\owedge \otimes \Id)\,\scrR^k
 &= (\nabla^k - \jet^k)\rmR
    + \frac{k+1}{k+3}(\owedge \otimes \Id)\,\scrR^k
    + \jet^k \rmR \\
 &= B_1(\nabla^{\le k-2} \rmR) + B_2(\nabla^{\le k-2} \rmR) ,
\end{align}
where \(B_2(\nabla^{\le k-2}\rmR)\) is the negative of
\eqref{eq:Young_Symmetrizer_modifiziert_1} divided by \(2(k+2)(k+3)\).
Hence we can set
\[
 B(\nabla^{\le k-2} \rmR) := B_1(\nabla^{\le k-2} \rmR) + B_2(\nabla^{\le k-2} \rmR)
\]
as claimed.
\end{proof}

The following example gives explicit formulas for
\(\nabla^{\le k-2} \rmR\) and at the same time shows the explicit
description of \(\nabla^{\le k}\rmR\) through \(\scrR^{\le k}\).

\bigskip
\begin{example}\label{ex:tilde_B}
\begin{enumerate}
\item[(a)]
For \(k = 2\) we obtain from
Proposition~\ref{p:Young_Symmetrizer_modifiziert_1}
\begin{equation*}
\begin{array}{l}
\big( \tfrac{1}{80}\,\rmS^\star_{\scaleto{
\begin{array}{|c|c|c|c|c|c|}
\cline{1-4}
1 & 3 & 5 & 6\\
\cline{1-4}
2 & 4 & \multicolumn{1}{c}{}\\
\cline{1-2}
\end{array}}{15pt}}\nabla^2_{X_5, X_6}
- \jet^2_{X_5, X_6} \big) \rmR_{X_1, X_2, X_3, X_4} \\
\quad= \tfrac{1}{80} \sum_{\sigma \in \rmS_{\{5, 6\}}} \Big( \tfrac{1}{2}\,
\rmS^\star_{\scaleto{
\begin{array}{|c|c|c|}
\cline{1-3}
1 & 3 & \sigma(5)\\
\cline{1-3}
2 & 4 & \multicolumn{1}{c}{}\\
\cline{1-2}
\end{array}}{15pt}}
\big( \rmR_{X_{\sigma(5)}, X_{\sigma(6)}}\,\rmR \big)_{X_1, X_2, X_3, X_4} \\
\quad\quad + 2\,\rmS_{\scaleto{
\begin{array}{|c|}
\cline{1-1}
1 \\
\cline{1-1}
2 \\
\cline{1-1}
\end{array}}{15pt}}
\big(\rmR_{X_1, X_{\sigma(5)}}\rmR \big)_{X_{\sigma(6)}, X_2, X_3, X_4} \\
\quad\quad - \rmS_{\scaleto{
\begin{array}{|c|}
\cline{1-1}
1 \\
\cline{1-1}
2 \\
\cline{1-1}
\end{array}}{15pt}}\,\rmS_{\scaleto{
\begin{array}{|c|}
\cline{1-1}
3 \\
\cline{1-1}
4 \\
\cline{1-1}
\end{array}}{15pt}}
\big(\rmR_{X_1, X_3}\rmR \big)_{X_{\sigma(5)}, X_2, X_{\sigma(6)}, X_4} \Big)
\end{array}
\end{equation*}
Clearly,
\begin{equation*}
\jet^2_{X_5, X_6} \rmR_{X_1, X_2, X_3, X_4}
= \nabla^2_{X_5, X_6} \rmR_{X_1, X_2, X_3, X_4}
- \tfrac{1}{2} \big(\rmR_{X_5, X_6} \rmR \big)_{X_1, X_2, X_3, X_4} .
\end{equation*}
Furthermore,
\begin{align*}
\rmS^\star_{\scaleto{
\begin{array}{|c|c|c|c|c|}
\cline{1-4}
1 & 3 & 5 & 6\\
\cline{1-4}
2 & 4 & \multicolumn{2}{c}{}\\
\cline{1-2}
\end{array}}{15pt}}\nabla^2_{X_5, X_6} \rmR_{X_1, X_2, X_3, X_4}
&= -48\,\rmS_{\scaleto{
\begin{array}{|c|}
\cline{1-1}
1 \\
\cline{1-1}
2 \\
\cline{1-1}
\end{array}}{15pt}}\,\rmS_{\scaleto{
\begin{array}{|c|}
\cline{1-1}
3 \\
\cline{1-1}
4 \\
\cline{1-1}
\end{array}}{15pt}}\, \scrR^2_{X_1, X_3, X_5, X_6; X_2, X_4} ,\\
\rmR_{X_1, X_2, X_3, X_4}
&= -\tfrac{1}{3}\,\rmS_{\scaleto{
\begin{array}{|c|}
\cline{1-1}
1 \\
\cline{1-1}
2 \\
\cline{1-1}
\end{array}}{15pt}}\,\rmS_{\scaleto{
\begin{array}{|c|}
\cline{1-1}
3 \\
\cline{1-1}
4 \\
\cline{1-1}
\end{array}}{15pt}}\, \scrR_{X_1, X_3; X_2, X_4} .
\end{align*}
According to \eqref{eq:YS1} and \eqref{eq:Fiedler}, it is clear how to
express \(\nabla^2 \rmR\) and \(\nabla^{\le 2}\rmR\) in terms of
\(\scrR^{\le 2}\).

\item[(b)]
For \(k = 3\), Proposition~\ref{p:Young_Symmetrizer_modifiziert_1}
implies that
\begin{align*}
&\big( \tfrac{1}{360}\,\rmS^\star_{\scaleto{
\begin{array}{|c|c|c|c|c|c|}
\cline{1-5}
1 & 3 & 5 & 6 & 7\\
\cline{1-5}
2 & 4 & \multicolumn{1}{c}{}\\
\cline{1-2}
\end{array}}{15pt}}\nabla^3_{X_5, X_6, X_7}
- \jet^3_{X_5, X_6, X_7} \big) \rmR_{X_1, X_2, X_3, X_4} \\
&= \tfrac{1}{60}\,\cyclic_{5, 6, 7}
\Big( \rmS^\star_{\scaleto{
\begin{array}{|c|c|c|}
\cline{1-3}
1 & 3 & 5\\
\cline{1-3}
2 & 4 & \multicolumn{1}{c}{}\\
\cline{1-2}
\end{array}}{15pt}}
\big( \jet^3_{X_5, X_6, X_7} - \jet^{2}_{X_6, X_7} \nabla_{X_5} \big)
\rmR_{X_1, X_2, X_3, X_4} \\
&\quad - \rmS^\star_{\scaleto{
\begin{array}{|c|c|}
\cline{1-2}
1 & 3 \\
\cline{1-2}
2 & 4 \\
\cline{1-2}
\end{array}}{15pt}}
\big( \nabla^3_{X_5, X_1, X_3} - \jet^3_{X_1, X_3, X_5} \big)
\rmR_{X_6, X_2, X_7, X_4} \Big) \\
&\quad + \sum_{\sigma \in \rmS_{\{5, 6, 7\}}} \Big(
2\,\rmS_{\scaleto{
\begin{array}{|c|}
\cline{1-1}
1 \\
\cline{1-1}
2 \\
\cline{1-1}
\end{array}}{15pt}}
\big( \nabla_{X_{\sigma(5)}} \rmR_{X_1, X_{\sigma(6)}} \rmR \big)_{X_{\sigma(7)}, X_2, X_3, X_4} \\
&\quad - \rmS_{\scaleto{
\begin{array}{|c|}
\cline{1-1}
1 \\
\cline{1-1}
2 \\
\cline{1-1}
\end{array}}{15pt}}\,\rmS_{\scaleto{
\begin{array}{|c|}
\cline{1-1}
3 \\
\cline{1-1}
4 \\
\cline{1-1}
\end{array}}{15pt}}
\big( \nabla_{X_{\sigma(5)}} \rmR_{X_1, X_3} \rmR \big)_{X_{\sigma(6)}, X_2, X_{\sigma(7)}, X_4}
\Big) .
\end{align*}
Here, for example,
\(\big( \jet^3_{X_5, X_6, X_7} - \jet^{2}_{X_6, X_7} \nabla_{X_5} \big)
 \rmR_{X_1, X_2, X_3, X_4}\)
is given by the negative of the right-hand side of~\eqref{eq:jet^3psi_1}
in Section~\ref{se:explicit_calculations}, where
\(\bbE := \scrC_0^\star(TM)\) is the vector bundle of algebraic
curvature tensors with the induced connection and \(\psi := \rmR\).
Similarly, the term
\(\big( \nabla^3_{X_5, X_1, X_3} - \jet^3_{X_1, X_3, X_5} \big)
 \rmR_{X_1, X_2, X_3, X_4}\)
corresponds to~\eqref{eq:jet^3psi_2}. Also recall that the Leibniz rule
must be applied to the terms \(\nabla_{X_i} \rmR_{X_j, X_k} \rmR\) as in
\eqref{eq:apply_Leibniz_0}.

\item[(c)]
For \(k = 4\), we proceed in a similar way and use~\eqref{eq:jet^4_4}
to obtain a description of \(\nabla^4 \rmR\) in terms of
\(\scrR^{\le 4}\) up to terms in \(\nabla^{\le 2} \rmR\). For an
explicit expression of \(\nabla^4 \rmR\) in terms of \(\scrR^{\le 4}\)
we have to turn to~(a).

\item[(d)]
Similarly, for \(k = 5\), use~\eqref{eq:jet^5_3} and~(b).
\end{enumerate}
\end{example}

\appendix
\section{Taylor expansion of the parallel transport}
\label{se:special_jet_formula}

Let a Riemannian manifold \((M, g)\) with Levi–Civita connection
\(\nabla\) and a curve \(c : \R \to M\) with \(c(0) = p\) be given. By
definition, the parallel transport \(\parallel_0^t(c) : T_p M \to
T_{c(t)} M\) is the fundamental solution of the ODE
\(\frac{\nabla}{\mathrm{d}t} Y(t) \equiv 0\). Equivalently,
\(\frac{\nabla}{\mathrm{d}t} \parallel_0^t(c) Y \equiv 0\) and
\(\parallel_0^0 Y = Y\) for all \(Y \in T_p M\), i.e.,~\(Y\) is
transported parallelly from \(p\) to \(c(t)\) along \(c\) for each \(t\).

The first goal of this section is to compare the simple jet formula
given in Lemma~\ref{le:Simple_Jet-Formula} with the special jet formula
from \cite[Ch.~3]{W1} in the version given in \cite[p.~32 (4.2)]{W1}:
\begin{equation}\label{eq:Special_Jet-Formula}
\begin{array}{ll}
\nabla^k_{X, \ldots, X, Y}\,\psi
&= \jet^{k}_{X, \ldots, X, Y}\,\psi
+ \sum_{r=2}^{k-1} \binom{k-1}{r}
\jet^{k-r}_{X, \ldots, X, \Phi_r(X) Y}\psi \\
&\quad + \sum_{r=1}^{k-1} \binom{k-1}{r}
\Omega^\mathbb{E}_r(X)_Y\,
\jet^{k-r-1}_{X, \ldots, X}\psi
\end{array}
\end{equation}
Here \(\Phi_r\) and \(\Omega^\mathbb{E}_r\) are tensors describing the
Levi–Civita connection \(\nabla\) and the linear connection
\(\nabla^\mathbb{E}\), respectively, in the “most natural” gauge
related to \(p\), namely with respect to geodesic normal coordinates and
the trivialization of \(\mathbb{E}\) obtained by using parallel
displacement along the radial geodesics emanating from \(p\).

More precisely, \(\frac{1}{r!} \Phi_r \in \Sym^r T^*M \otimes
\End(T_p M)\) is by definition the \(r\)th coefficient of the Taylor
expansion
\[
\Phi(X) Y \underset{X \to 0}{\sim}
\sum_{r=0}^\infty \frac{1}{r!} \Phi_r(X) Y
\]
of the parallel transport map
\begin{equation}\label{eq:Phi}
  \Phi : U \to \End(T_p M), \quad
  X \mapsto \Phi(X) : T_p M \xrightarrow{\parallel_0^1 \gamma_X}
  T_{\exp^M_p(X)} M \xrightarrow{(D_X \exp^M_p)^{-1}}
  T_X T_p M \cong T_p M
\end{equation}
Here \(U\) is some open star-shaped neighborhood of \(0\) in \(T_p M\)
where the exponential map \(\exp_p^M\) defines an isomorphism onto
\(\exp_p^M(U)\), \(\parallel_0^1 \gamma_X\) is the parallel transport in
\(TM\) from \(0\) to \(1\) along the geodesic \(\gamma_X(t) = \exp^M_p(t
X)\) emanating from \(p\), and \((D_X \exp^M_p)^{-1}\) is the inverse of
the differential of the exponential map \(\exp^M_p : T_p M \to M\).

Similarly, we have a Taylor expansion
\[
\Omega^\mathbb{E}(X)_Y \psi \underset{X \to 0}{\sim}
\sum_{r=0}^\infty \frac{1}{r!} \Omega^\mathbb{E}_r(X)_Y \psi
\]
where \(\Omega^\mathbb{E}(X)_Y \psi := \omega^\mathbb{E}(X)_{\Phi(X) Y}
\psi\) and \(\omega^\mathbb{E} : U \to T^*_p M \otimes
\End(\mathbb{E}_p)\), \(X \mapsto \omega^\mathbb{E}(X)\), is the 1-form
describing the linear connection \(\nabla^\mathbb{E}\) via
\[
\nabla_Y \psi(X) =
\frac{\partial}{\partial Y} \psi(X) + \omega^\mathbb{E}(X)_Y \psi
\]
with respect to the local trivialization
\[
U \times \mathbb{E}_p \overset{\sim}{\longrightarrow}
\mathbb{E}|_{\exp^M_p(U)}, \quad
(X, \psi) \mapsto (\parallel_0^1 \gamma_X)^\mathbb{E} \psi
\]
obtained by parallel translation \((\parallel_0^1 \gamma_X)^\mathbb{E}\)
of \(\mathbb{E}_p\) along the radial geodesics \(\gamma_X\) from \(p\).

Comparing \eqref{eq:Simple_Jet-Formula} and
\eqref{eq:Special_Jet-Formula},
\begin{equation}\label{eq:Simple_vs_Special}
\begin{array}{ll}
\sum_{j=1}^{k-1} \frac{j}{k}\,
\nabla^{j-1}_{X, \ldots, X} \rmR_{X, Y}\,
\nabla^{k-j-1}_{X, \ldots, X} \psi
&= \sum_{r=2}^{k-1} \binom{k-1}{r}
\jet^{k-r}_{X, \ldots, X, \Phi_r(X) Y} \psi \\
&\quad + \sum_{r=1}^{k-1} \binom{k-1}{r}
\Omega^\mathbb{E}_r(X)_Y\,
\jet^{k-r-1}_{X, \ldots, X} \psi
\end{array}
\end{equation}
From this we can obtain the Taylor expansion of both \(\Phi\) and
\(\Omega^\mathbb{E}\) up to arbitrary order by using the Leibniz rule
for the Ricci identity \eqref{eq:Ricci_refined}. For example,
\begin{align*}
\Phi(X) Y
&= Y - \frac{1}{6} \rmR_{X, Y} X
- \frac{1}{12} (\nabla_X \rmR)_{X, Y} X
- \frac{1}{40} (\nabla^2_{X, X} \rmR)_{X, Y} X
- \frac{7}{360} (\rmR_{\rmR_{X, Y} X, X} X), \\
\Phi^{-1}(X) Y
&= Y + \frac{1}{6} \rmR_{X, Y} X
+ \frac{1}{12} (\nabla_X \rmR)_{X, Y} X
+ \frac{1}{40} (\nabla^2_{X, X} \rmR)_{X, Y} X
- \frac{1}{120} (\rmR_{\rmR_{X, Y} X, X} X), \\
\Omega^\mathbb{E}(X)_Y
&= \frac{1}{2} \rmR^\mathbb{E}_{X, Y}
+ \frac{1}{3} (\nabla_X \rmR^\mathbb{E})_{X, Y}
+ \frac{1}{8} (\nabla^2_{X, X} \rmR^\mathbb{E})_{X, Y}
+ \frac{1}{24} \rmR^\mathbb{E}_{\rmR_{X, Y} X, X} \\
&\quad + \frac{1}{30} (\nabla^3_{X, X, X} \rmR^\mathbb{E})_{X, Y}
+ \frac{1}{45} (\nabla_X \rmR^\mathbb{E})_{\rmR_{X, Y} X, X}
+ \frac{1}{40} \rmR^\mathbb{E}_{(\nabla_X \rmR)_{X, Y} X, X}, \\
\omega^\mathbb{E}(X)_Y
&= \frac{1}{2} \rmR^\mathbb{E}_{X, Y}
+ \frac{1}{3} (\nabla_X \rmR^\mathbb{E})_{X, Y}
+ \frac{1}{8} (\nabla^2_{X, X} \rmR^\mathbb{E})_{X, Y}
- \frac{1}{24} \rmR^\mathbb{E}_{\rmR_{X, Y} X, X} \\
&\quad + \frac{1}{30} (\nabla^3_{X, X, X} \rmR^\mathbb{E})_{X, Y}
- \frac{1}{30} (\nabla_X \rmR^\mathbb{E})_{\rmR_{X, Y} X, X}
- \frac{1}{60} \rmR^\mathbb{E}_{(\nabla_X \rmR)_{X, Y} X, X}
\end{align*}
These are the Taylor polynomials of \(\Phi\), \(\Phi^{-1}\),
\(\Omega^\mathbb{E}\), and \(\omega^\mathbb{E}\) of order four.

Moreover, in Section~\ref{se:Taylor_expansion} we find an explicit
formula for the Taylor coefficients of \(\Phi^{-1}\). Since \(\Phi(X)\)
is by definition the inverse of \(\Phi^{-1}(X)\) for each \(X \in U\),
the Taylor series of \(\Phi\) can alternatively be obtained from that of
\(\Phi^{-1}\) by formal inversion (for example, the equality
\(\frac{1}{6} \cdot \frac{1}{6} - \frac{1}{120} = \frac{7}{360}\)
relates the terms quadratic in \(\scrR\)). However, a similarly simple
formula for the coefficients of \(\Omega^\mathbb{E}\) and
\(\omega^\mathbb{E}\) is seemingly not known.

\subsection{Explicit calculations for the jet formula from
  Proposition \ref{p:Simple_Jet-Formula}}
\label{se:explicit_calculations}

In the following example, we write out \eqref{eq:jet^k-nabla^k} in
detail for small values of \(k\) and, as a byproduct, determine the
coefficients of the Taylor expansions of \(\Phi\) and \(\Omega^\bbE\).

\bigskip
\begin{example}\label{ex:k=2345}
\begin{enumerate}

\item For \(k=2\), it is immediate that
\begin{equation}\label{eq:special_jet_k=2}
  \nabla^2_{X,Y}\psi - \jet^{2}_{X,Y}\psi
  = \frac{1}{2}\rmR^\bbE_{X,Y}\psi
\end{equation}
Hence, by \eqref{eq:Simple_vs_Special},
\begin{equation}\label{eq:Omega_1}
  \Omega^\bbE_1(X)_Y = \frac{1}{2}\rmR^\bbE_{X,Y}
\end{equation}

\item For \(k=3\), \eqref{eq:Simple_Jet-Formula} gives
\begin{equation}\label{eq:special_jet_k=3}
  \nabla^{3}_{X,X,Y}\psi - \jet^{3}_{X,X,Y}\psi
  = \frac{1}{3}\rmR_{X,Y}\nabla_{X}\psi
    + \frac{2}{3}(\nabla_X \rmR)_{X,Y}\psi .
\end{equation}
Here we have not yet applied the Leibniz rule. According to
\eqref{eq:Ricci_refined},
\begin{equation}\label{eq:apply_Leibniz_0}
  \nabla_{X}\big(\rmR_{X,Y}\psi\big)
  = (\nabla_{X}\rmR^\bbE)_{X,Y}\psi
    + \rmR^\bbE_{X,Y}\nabla_X\psi .
\end{equation}
Hence,
\begin{equation}\label{eq:apply_Leibniz_1}
  \rmR_{X,Y}\nabla_{X}\psi
  + 2\nabla_{X}\rmR_{X,Y}\psi
  = 2(\nabla_{X}\rmR^\bbE)_{X,Y}\psi
    + 3\rmR^\bbE_{X,Y}\nabla_{X}\psi
    - \nabla_{\rmR_{X,Y}X}\psi .
\end{equation}
Substituting into \eqref{eq:special_jet_k=3} yields
\begin{equation}\label{eq:jet^3psi_1}
  \nabla^3_{X,X,Y}\psi
  = \jet^{3}_{X,X,Y}\psi
    + \frac{2}{3}(\nabla_{X}\rmR^\bbE)_{X,Y}\psi
    + \rmR^\bbE_{X,Y}\nabla_{X}\psi
    - \frac{1}{3}\nabla_{\rmR_{X,Y}X}\psi .
\end{equation}

After polarization we obtain
\begin{equation}\label{eq:jet^3psi_2}
\begin{aligned}
  &\jet^2_{X_1,X_2}\nabla_{X_3}\psi
   - \jet^{3}_{X_1,X_2,X_3}\psi\\
  &= \frac{1}{2}\sum_{\sigma\in \rmS_2}
  \Big(
    \frac{2}{3}(\nabla_{X_{\sigma(1)}} \rmR^\bbE)_{X_{\sigma(2)},X_3}
    + \rmR^\bbE_{X_{\sigma(1)},X_3}\nabla_{X_{\sigma(2)}}
    - \frac{1}{3}\nabla_{\rmR_{X_{\sigma(1)},X_3}X_{\sigma(2)}}
  \Big)\psi 
\end{aligned}
\end{equation}
Comparing with \eqref{eq:Simple_vs_Special} gives
\begin{equation}\label{eq:Phi_2}
  \Phi_2(X)Y = -\frac{1}{3}\rmR_{X,Y}X
\qquad
  \Omega_2^\bbE(X)_Y = \frac{2}{3}(\nabla_X \rmR^\bbE)_{X,Y}
\end{equation}
Also,
\begin{equation}\label{eq:apply_Leibniz_3}
 \rmR_{X,Y}\nabla_{Z}\psi
 =
 \rmR^\bbE_{X,Y}\nabla_{Z}\psi
 - \nabla_{\rmR_{X,Y}Z}\psi
\end{equation}
Therefore, applying \eqref{eq:special_jet_k=2} to the vector bundle
\(T^*M \otimes \bbE\) and the section \(\tilde\psi := \nabla\psi\), and
using \eqref{eq:jet^3psi_1}, we obtain
\begin{align}\notag
 &\nabla^3_{X_1,X_2,X_3}\psi
  - \jet^{3}_{X_1,X_2,X_3}\psi \\[-1ex]\notag
 &\quad=
 \big(\jet^2_{X_1,X_2}\nabla_{X_3}
      - \jet^3_{X_1,X_2,X_3}\big)\psi
 + \big(\nabla^3_{X_1,X_2,X_3}
      - \jet^2_{X_1,X_2}\nabla_{X_3}\big)\psi \\[0.5ex]
 &\quad= \frac{1}{2} \sum_{\sigma \in \rmS_2}
 \Big(\frac{2}{3}(\nabla_{X_{\sigma(1)}}\rmR^\bbE)_{X_{\sigma(2)},X_3}
 + \rmR^\bbE_{X_{\sigma(1)},X_3}\nabla_{X_{\sigma(2)}}
 - \frac{1}{3}\nabla_{\rmR_{X_{\sigma(1)},X_3}X_{\sigma(2)}}\Big)\psi
 \\[-1ex]\notag
 &\qquad\quad
 + \frac{1}{2}\big(\rmR^\bbE_{X_1,X_2}\nabla_{X_3}
 - \nabla_{\rmR_{X_1,X_2}X_3}\big)\psi
\end{align}
cf. (4.3) from \cite[p.~35]{W1}.

\item For \(k = 4\), \eqref{eq:Simple_Jet-Formula} gives
\begin{equation}\label{eq:jet^4_1}
 \nabla^{4}_{X,X,X,Y}\psi - \jet^{4}_{X,X,X,Y}\psi
 =
 \frac{1}{4} \big(\rmR_{X,Y}\nabla^2_{X,X}
 + 2\nabla_{X}\rmR_{X,Y}\nabla_{X}
 + 3\nabla^2_{X,X}\rmR_{X,Y}\big)\psi
\end{equation}
Moreover, incorporating the Leibniz rule as in
\eqref{eq:Ricci_refined},
\begin{align*}
 \rmR_{X,Y}\nabla^2_{X,X}\psi
 &= \big(\rmR^\bbE_{X,Y}\nabla^2_{X,X}
 - \nabla^2_{\rmR_{X,Y}X,X}
 - \nabla^2_{X,\rmR_{X,Y}X}\big)\psi \\
 \nabla_{X}\rmR_{X,Y}\nabla_{X}\psi
 &= \big((\nabla_{X}\rmR^\bbE)_{X,Y}\nabla_{X}
 - \nabla_{(\nabla_{X}\rmR)_{X,Y}X}
 + \rmR^\bbE_{X,Y}\nabla^2_{X,X}
 - \nabla^2_{X,\rmR_{X,Y}X}\big)\psi \\
 \nabla^2_{X,X}\rmR_{X,Y}\psi
 &= \big((\nabla^2_{X,X}\rmR^\bbE)_{X,Y}
 + 2(\nabla_{X}\rmR^\bbE)_{X,Y}\nabla_{X}
 + \rmR^\bbE_{X,Y}\nabla^2_{X,X}\big)\psi
\end{align*}
We conclude that
\begin{align}\label{eq:jet^4_2}
 \nabla^{4}_{X,X,X,Y}\psi - \jet^{4}_{X,X,X,Y}\psi
 &= \frac{1}{4} \big(3(\nabla^2_{X,X}\rmR^\bbE)_{X,Y}
 + 8(\nabla_{X}\rmR^\bbE)_{X,Y}\nabla_{X}
 + 6\rmR^\bbE_{X,Y}\nabla^2_{X,X} \notag \\
 &\quad - 4\jet^2_{\rmR_{X,Y}X,X}
 - \rmR^\bbE_{X,\rmR_{X,Y}X}
 - 2\nabla_{(\nabla_X\rmR)_{X,Y}X}\big)\psi
\end{align}
Now, \eqref{eq:Simple_vs_Special} gives
\begin{align}
 \Phi_3(X)Y
 &= -\frac{1}{2}(\nabla_X\rmR)_{X,Y}X , \label{eq:Phi_3} \\
 \Omega_3^\bbE(X)_Y\psi
 &= \frac{3}{4}(\nabla^2_{X,X}\rmR^\bbE)_{X,Y}\psi
 + \frac{1}{4}\rmR^\bbE_{\rmR_{X,Y}X,X}\psi \label{eq:Omega_3}
\end{align}
By polarization in \(X\), \eqref{eq:jet^4_2} becomes
\begin{align}\label{eq:jet^4_3}
 &\jet^3_{X_1,X_2,X_3}\nabla_{X_4}\psi
 - \jet^{4}_{X_1,X_2,X_3,X_4}\psi \notag \\
 &\quad= \frac{1}{6} \sum_{\sigma \in \rmS_3}
 \Big(\frac{3}{4}(\nabla^2_{X_{\sigma(1)},X_{\sigma(2)}}\rmR^\bbE)_{X_{\sigma(3)},X_4}
 + 2(\nabla_{X_{\sigma(1)}}\rmR^\bbE)_{X_{\sigma(2)},X_4}\nabla_{X_{\sigma(3)}} \notag \\
 &\qquad\quad
 + \frac{3}{2}\rmR^\bbE_{X_{\sigma(1)},X_4}\nabla^2_{X_{\sigma(2)},X_{\sigma(3)}}
 - \jet^2_{\rmR_{X_{\sigma(1)},X_4}X_{\sigma(2)},X_{\sigma(3)}} \notag \\
 &\qquad\quad
 - \frac{1}{4}\rmR^\bbE_{X_{\sigma(1)},\rmR_{X_{\sigma(2)},X_4}X_{\sigma(3)}}
 - \frac{1}{2}\nabla_{(\nabla_{X_{\sigma(1)}}\rmR)_{X_{\sigma(2)},X_4}X_{\sigma(3)}}\Big)\psi
\end{align}
Following the proof of Proposition~\ref{p:Simple_Jet-Formula}, we
finally obtain
\begin{align}\notag
 &\big( \nabla^4_{X_1, X_2, X_3, X_4}
       - \jet^{4}_{X_1, X_2, X_3, X_4} \big) \psi\\
 &\begin{array}{llll}
   =\; &\big( \jet^3_{X_1, X_2, X_3} \nabla_{X_4}
           - \jet^4_{X_1, X_2, X_3, X_4} \big) \psi\\
       &+ \big( \jet^2_{X_1, X_2} \nabla^2_{X_3, X_4}
           - \jet^3_{X_1, X_2, X_3} \nabla_{X_4} \big) \psi\\
       &+ \big( \nabla^4_{X_1, X_2, X_3, X_4}
           - \jet^2_{X_1, X_2} \nabla^2_{X_3, X_4} \big) \psi
  \end{array}\\\label{eq:jet^4_4}
 &\begin{array}{ll}
   =\; \frac{1}{6} \sum_{\sigma \in \rmS_3}
   &\big( \tfrac{3}{4} (\nabla^2_{X_{\sigma(1)}, X_{\sigma(2)}}
        \rmR^\bbE)_{X_{\sigma(3)}, X_4}
        + 2 (\nabla_{X_{\sigma(1)}} \rmR^\bbE)_{X_{\sigma(2)}, X_4}
          \nabla_{X_{\sigma(3)}} \\
   &\quad + \tfrac{3}{2} \rmR^\bbE_{X_{\sigma(1)}, X_4}
          \nabla^2_{X_{\sigma(2)}, X_{\sigma(3)}}
        + \tfrac{1}{4} \rmR^\bbE_{\rmR_{X_{\sigma(1)}, X_4}
          X_{\sigma(2)}, X_{\sigma(3)}} \\
   &\quad - \jet^2_{\rmR_{X_{\sigma(1)}, X_4} X_{\sigma(2)}, X_{\sigma(3)}}
        - \tfrac{1}{2} \nabla_{(\nabla_{X_{\sigma(1)}} \rmR)_{X_{\sigma(2)},
          X_4} X_{\sigma(3)}} \big) \psi\\
   &+ \frac{1}{2} \sum_{\sigma \in \rmS_2}
     \big( \tfrac{2}{3} (\nabla_{X_{\sigma(1)}} \rmR^\bbE)_{X_{\sigma(2)},
          X_3} \nabla_{X_4}
        + \rmR^\bbE_{X_{\sigma(1)}, X_3}
          \nabla^2_{X_{\sigma(2)}, X_4} \\
   &\quad - \tfrac{1}{3} \nabla^2_{\rmR_{X_{\sigma(1)}, X_3} X_{\sigma(2)},
          X_4}
        - \tfrac{2}{3} \nabla_{(\nabla_{X_{\sigma(1)}} \rmR)_{X_{\sigma(2)},
          X_3} X_4} \\
   &\quad - \nabla^2_{X_{\sigma(1)}, \rmR_{X_{\sigma(2)}, X_3} X_4} \big) \psi\\
   &+ \frac{1}{2}
     \big( \rmR^\bbE_{X_1, X_2} \nabla^2_{X_3, X_4}
        - \nabla^2_{\rmR_{X_1, X_2} X_3, X_4}
        - \nabla^2_{X_3, \rmR_{X_1, X_2} X_4} \big) \psi
  \end{array}
\end{align}
By rewriting in the above formula all second-order covariant
derivatives \(\nabla^2\) that act directly on \(\psi\) as
\(\jet^2 + \tfrac{1}{2} \rmR^\bbE\) (e.g.,
\(\rmR^\bbE_{X_1,X_2} \nabla^2_{X_3,X_4} \psi
= \rmR^\bbE_{X_1,X_2} \jet^2_{X_3,X_4} \psi
+ \tfrac{1}{2} \rmR^\bbE_{X_1,X_2} \rmR^\bbE_{X_3,X_4} \psi\)) —
except for the term \(\tfrac{1}{4} \sum_{\sigma \in \rmS_3}
\rmR^\bbE_{X_{\sigma(1)},X_4} \nabla^2_{X_{\sigma(2)},X_{\sigma(3)}}
\psi\) — it is straightforward to check that \eqref{eq:jet^4_4} is, in
fact, consistent with the expression of \(\nabla^4 \psi\) in
\(\jet^{\le 4} \psi\) obtained by summing up the terms related to the
coefficients \eqref{eq:Phi_2}, \eqref{eq:Phi_3} and \eqref{eq:Omega_1},
\eqref{eq:Omega_3} of the Taylor polynomials of order three of \(\Phi\)
and \(\Omega\), respectively, via the thirty jet forests of order four
with feedback as described in \cite[Lemma~4.2]{W1}.

\item For \(k=5\), from \eqref{eq:Simple_Jet-Formula} we obtain
\begin{equation}\label{eq:jet^5_1}
\begin{aligned}
\big(\nabla^{5}_{X,X,X,X,Y} - \jet^{5}_{X,X,X,X,Y}\big)\psi
= \frac{1}{5}\Big(
  \rmR_{X,Y}\nabla^3_{X,X,X}
 + 2\nabla_{X}\rmR_{X,Y}\nabla^2_{X,X} \\
 + 3\nabla^2_{X,X}\rmR_{X,Y}\nabla_X
 + 4\nabla^3_{X,X,X}\rmR_{X,Y}
\Big)\psi 
\end{aligned}
\end{equation}
Incorporating the Leibniz rule as in \eqref{eq:Ricci_refined} yields
\begin{align*}
\rmR_{X,Y}\nabla^3_{X,X,X}\psi
&= \big(
    \rmR^\bbE_{X,Y}\nabla^3_{X,X,X}
  - \nabla^3_{\rmR_{X,Y}X,X,X}
  - \nabla^3_{X,\rmR_{X,Y}X,X}
  - \nabla^3_{X,X,\rmR_{X,Y}X}
  \big)\psi \\
\nabla_{X}\rmR_{X,Y}\nabla^2_{X,X}\psi
&= \big(
    (\nabla_X \rmR^\bbE)_{X,Y}\nabla^2_{X,X}
  - \nabla^2_{(\nabla_X \rmR)_{X,Y}X,X}
  - \nabla^2_{X,(\nabla_X \rmR)_{X,Y}X} \\
&\qquad\quad
  + \rmR^\bbE_{X,Y}\nabla^3_{X,X,X}
  - \nabla^3_{X,\rmR_{X,Y}X,X}
  - \nabla^3_{X,X,\rmR_{X,Y}X}
  \big)\psi \\
\nabla^2_{X,X}\rmR_{X,Y}\nabla_X\psi
&= \big(
    (\nabla^2_{X,X} \rmR^\bbE)_{X,Y}\nabla_X
  + 2(\nabla_X \rmR^\bbE)_{X,Y}\nabla^2_{X,X}
  + \rmR^\bbE_{X,Y}\nabla^3_{X,X,X} \\
&\qquad\quad
  - \nabla_{(\nabla^2_{X,X} \rmR)_{X,Y}X}
  - 2\nabla^2_{X,(\nabla_X \rmR)_{X,Y}X}
  - \nabla^3_{X,X,\rmR_{X,Y}X}
  \big)\psi \\
\nabla^3_{X,X,X}\rmR_{X,Y}\psi
&= \big(
    (\nabla^3_{X,X,X} \rmR^\bbE)_{X,Y}
  + 3(\nabla^2_{X,X} \rmR^\bbE)_{X,Y}\nabla_X
  + 3(\nabla_X \rmR^\bbE)_{X,Y}\nabla^2_{X,X} \\
&\qquad\quad
  + \rmR^\bbE_{X,Y}\nabla^3_{X,X,X}
  \big)\psi 
\end{align*}
Hence,
\begin{equation}\label{eq:jet^5_2}
\begin{aligned}
\big(\nabla^{5}_{X,X,X,X,Y} - \jet^{5}_{X,X,X,X,Y}\big)\psi
= \Big(
  \tfrac{4}{5}(\nabla^3_{X,X,X} \rmR^\bbE)_{X,Y}
 + 3(\nabla^2_{X,X} \rmR^\bbE)_{X,Y}\nabla_X \\
 + 4(\nabla_X \rmR^\bbE)_{X,Y}\nabla^2_{X,X}
 + 2\rmR^\bbE_{X,Y}\nabla^3_{X,X,X} \\
 - 2\jet^3_{\rmR_{X,Y}X,X,X}
 + \rmR^\bbE_{\rmR_{X,Y}X,X}\nabla_X \\
 - 2\jet^2_{(\nabla_X \rmR)_{X,Y}X,X}
 + \tfrac{3}{5}\rmR^\bbE_{(\nabla_X \rmR)_{X,Y}X,X}
 - \tfrac{3}{5}\nabla_{(\nabla^2_{X,X}\rmR)_{X,Y}X} \\
 + \tfrac{8}{15}(\nabla_X \rmR^\bbE)_{\rmR_{X,Y}X,X}
 - \tfrac{7}{15}\nabla_{\rmR_{\rmR_{X,Y}X,X}X}
\Big)\psi
\end{aligned}
\end{equation}
Therefore, \eqref{eq:Simple_vs_Special} gives
\begin{align}
&\Phi_4(X)Y = -\frac{3}{5}(\nabla^2_{X,X}\rmR)_{X,Y}X
               - \frac{7}{15}\big(\rmR_{\rmR_{X,Y}X,X}X\big)\\
&\Omega_4^\bbE(X)_Y
 = \frac{4}{5}(\nabla^3_{X,X,X}\rmR^\bbE)_{X,Y}
 + \frac{8}{15}(\nabla_X \rmR^\bbE)_{\rmR_{X,Y}X,X}
 + \frac{3}{5}\rmR^\bbE_{(\nabla_X \rmR)_{X,Y}X,X}
\end{align}

By polarization in \(X\), we obtain
\begin{equation}\label{eq:jet^5_3}
\begin{aligned}
&\big(\nabla^5_{X_1,X_2,X_3,X_4,X_5}
   - \jet^5_{X_1,X_2,X_3,X_4,X_5}\big)\psi \\
&= \frac{1}{24}\sum_{\sigma \in \rmS_4} \Big(
    \frac{4}{5}(\nabla^3_{X_{\sigma(1)},X_{\sigma(2)},X_{\sigma(3)}}
      \rmR^\bbE)_{X_{\sigma(4)},X_5}
  + 3(\nabla^2_{X_{\sigma(1)},X_{\sigma(2)}} \rmR^\bbE)_{X_{\sigma(3)},X_5}
      \nabla_{X_{\sigma(4)}} \\
&\qquad\qquad\qquad\qquad\quad
  + 4(\nabla_{X_{\sigma(1)}} \rmR^\bbE)_{X_{\sigma(2)},X_5}
      \nabla^2_{X_{\sigma(3)},X_{\sigma(4)}}
  + 2\rmR^\bbE_{X_{\sigma(1)},X_5}
      \nabla^3_{X_{\sigma(2)},X_{\sigma(3)},X_{\sigma(4)}} \\
&\qquad\qquad\qquad\qquad\quad
  - 2\jet^3_{\rmR_{X_{\sigma(1)},X_5}X_{\sigma(2)},
      X_{\sigma(3)},X_{\sigma(4)}}
  + \rmR^\bbE_{\rmR_{X_{\sigma(1)},X_5}X_{\sigma(2)},
      X_{\sigma(3)}}\nabla_{X_{\sigma(4)}} \\
&\qquad\qquad\qquad\qquad\quad
  - 2\jet^2_{(\nabla_{X_{\sigma(1)}} \rmR)_{X_{\sigma(2)},X_5}
      X_{\sigma(3)},X_{\sigma(4)}}
  + \frac{3}{5}\rmR^\bbE_{(\nabla_{X_{\sigma(1)}} \rmR)_{X_{\sigma(2)},X_5}
      X_{\sigma(3)},X_{\sigma(4)}} \\
&\qquad\qquad\qquad\qquad\quad
  - \frac{3}{5}\nabla_{(\nabla^2_{X_{\sigma(1)},X_{\sigma(2)}} \rmR)
      _{X_{\sigma(3)},X_5} X_{\sigma(4)}}
  + \frac{8}{15}(\nabla_{X_{\sigma(1)}} \rmR^\bbE)
      _{\rmR_{X_{\sigma(2)},X_5}X_{\sigma(3)},X_{\sigma(4)}} \\
&\qquad\qquad\qquad\qquad\quad
  - \frac{7}{15}\nabla_{\rmR_{\rmR_{X_{\sigma(1)},X_5}X_{\sigma(2)},
      X_{\sigma(3)}} X_{\sigma(4)}}
\Big)\psi \\
&\quad
+ \frac{1}{6}\sum_{\sigma \in \rmS_3} \Big(
    \frac{3}{4}(\nabla^2_{X_{\sigma(1)},X_{\sigma(2)}} \rmR^\bbE)
      _{X_{\sigma(3)},X_4}\nabla_{X_5}
  + 2(\nabla_{X_{\sigma(1)}} \rmR^\bbE)
      _{X_{\sigma(2)},X_4}\nabla^2_{X_{\sigma(3)},X_5} \\
&\qquad\qquad\qquad\qquad\quad
  + \frac{3}{2}\rmR^\bbE_{X_{\sigma(1)},X_4}
      \nabla^3_{X_{\sigma(2)},X_{\sigma(3)},X_5}
  - \frac{1}{4}\rmR^\bbE_{X_{\sigma(1)},\rmR_{X_{\sigma(2)},X_4}
      X_{\sigma(3)}}\nabla_{X_5} \\
&\qquad\qquad\qquad\qquad\quad
  - \frac{1}{2}\nabla^3_{\rmR_{X_{\sigma(1)},X_4}X_{\sigma(2)},
      X_{\sigma(3)},X_5}
  - \frac{1}{2}\nabla^3_{X_{\sigma(1)},\rmR_{X_{\sigma(2)},X_4}
      X_{\sigma(3)},X_5} \\
&\qquad\qquad\qquad\qquad\quad
  - \frac{1}{2}\nabla^2_{(\nabla_{X_{\sigma(1)}} \rmR)
      _{X_{\sigma(2)},X_4}X_{\sigma(3)},X_5}
  - \frac{3}{4}\nabla_{(\nabla^2_{X_{\sigma(1)},X_{\sigma(2)}} \rmR)
      _{X_{\sigma(3)},X_4} X_5} \\
&\qquad\qquad\qquad\qquad\quad
  - 2\nabla^2_{X_{\sigma(1)},(\nabla_{X_{\sigma(2)}} \rmR)
      _{X_{\sigma(3)},X_4} X_5}
  - \frac{3}{2}\nabla^3_{X_{\sigma(1)},X_{\sigma(2)},\rmR_{X_{\sigma(3)},
      X_4} X_5}\\
&\qquad\qquad\qquad\qquad\quad
  + \frac{1}{4}\nabla_{\rmR_{X_{\sigma(1)},\rmR_{X_{\sigma(2)},X_4}
      X_{\sigma(3)}} X_5}
\Big)\psi \\
&\quad
+ \frac{1}{2}\sum_{\sigma \in \rmS_2} \Big(
    \frac{2}{3}(\nabla_{X_{\sigma(1)}} \rmR^\bbE)
      _{X_{\sigma(2)},X_3}\nabla^2_{X_4,X_5}
  + \rmR^\bbE_{X_{\sigma(1)},X_3}
      \nabla^3_{X_{\sigma(2)},X_4,X_5} \\
&\qquad\qquad\qquad\qquad\quad
  - \frac{1}{3}\nabla^3_{\rmR_{X_{\sigma(1)},X_3} X_{\sigma(2)},
      X_4,X_5}
  - \nabla^3_{X_{\sigma(1)},X_4,\rmR_{X_{\sigma(2)},X_3} X_5}
  - \nabla^3_{X_{\sigma(1)},\rmR_{X_{\sigma(2)},X_3} X_4,X_5} \\
&\qquad\qquad\qquad\qquad\quad
  - \frac{2}{3}\nabla^2_{X_4,(\nabla_{X_{\sigma(1)}} \rmR)
      _{X_{\sigma(2)},X_3} X_5}
  - \frac{2}{3}\nabla^2_{(\nabla_{X_{\sigma(1)}} \rmR)
      _{X_{\sigma(2)},X_3} X_4,X_5}
\Big)\psi \\
&\quad
+ \frac{1}{2}\Big(
    \rmR^\bbE_{X_1,X_2}\nabla^3_{X_3,X_4,X_5}
  - \nabla^3_{\rmR_{X_1,X_2}X_3,X_4,X_5}
  - \nabla^3_{X_3,\rmR_{X_1,X_2}X_4,X_5}
  - \frac{1}{2}\nabla^3_{X_3,X_4,\rmR_{X_1,X_2}X_5}
\Big)\psi
\end{aligned}
\end{equation}

\end{enumerate}
\end{example}
\section{Taylor expansion of the metric in normal coordinates}
\label{se:Taylor_expansion}

To clarify the notion of a noncommutative polynomial in
Theorem~\ref{th:taylor_series}, consider the unital associative
\(\bbR\)-algebra
\[
  \scrA_{\text{univ}} := \bbR\!\left\langle \scrR^0, \scrR^1, \scrR^2,
  \ldots \right\rangle
\]
freely generated by a countable family \(\{\scrR^i\}_{i\ge0}\). It is
characterized by the universal property: for any unital associative
\(\bbR\)-algebra \(\scrA\) and any sequence
\((\tilde\scrR^i)_{i\ge0} \subset \scrA\), there exists a unique
homomorphism
\[
  \ev_{\tilde\scrR}\colon \scrA_{\text{univ}} \longrightarrow \scrA
\]
such that \(\ev_{\tilde\scrR}(\scrR^i)=\tilde\scrR^i\) for all \(i\).

Elements of \(\scrA_{\text{univ}}\) are finite \(\bbR\)-linear
combinations of words \(\scrR^I := \scrR^{i_1}\cdots \scrR^{i_r}\)
with \(I=(i_1,\ldots,i_r)\) and \(r\ge0\) (the empty word for \(r=0\) is
the unit). Evaluation is substitution:
\(\ev_{\tilde\scrR}(\scrR^I)=\tilde\scrR^{i_1}\cdots \tilde\scrR^{i_r}\).

Since the symmetrized \(k\)th covariant derivative
\(\scrR^k|_p\) of the curvature tensor of a Riemannian manifold
\((M,g)\) is a polynomial of degree \(k+2\) on \(T_pM\) with values in
\(\End{T_pM}\) (see \eqref{eq:k-ter_Jacobi}), we equip
\(\scrA_{\text{univ}}\) with the grading
\[
  \deg(\scrR^k)=k+2, \qquad
  \deg(\scrR^I)=i_1+\cdots+i_r+2r
\]
and call an expression \emph{homogeneous} if it is supported in a
single total degree.

Next, let \(V\) be a vector space and set
\(\scrA := \Sym^\bullet V^* \otimes \End{V}\), graded by the
polynomial degree on \(\Sym^\bullet V^*\) and with multiplication
\[
  (h_1 \otimes a_1)\cdot(h_2 \otimes a_2) := (h_1 h_2) \otimes
  (a_1 \circ a_2), \qquad \deg(h \otimes a) := \deg(h)
\]
for \(h,h_1,h_2 \in \Sym^\bullet V^*\) and \(a,a_1,a_2 \in \End{V}\).
Given \(\tilde\scrR^i \in \Sym^{i+2}V^* \otimes \End{V}\) for \(i\ge0\),
any \(Q \in \scrA_{\text{univ}}\) evaluates to
\[
  Q(\tilde\scrR^0,\tilde\scrR^1,\ldots) \in \Sym^\bullet V^* \otimes
  \End{V}
\]
and for \(X \in V\),
\[
  Q(\tilde\scrR^0,\tilde\scrR^1,\ldots)(X) =
  Q\bigl(\tilde\scrR^0(X),\tilde\scrR^1(X),\ldots\bigr)
\]
which is a polynomial in \(X\) of the same total degree as \(Q\).

\subsection{Taylor expansion of the backward parallel transport}
\label{se:backward}

For a smooth curve \(c\), write
\(\parallel_{s}^{t}(c)\colon T_{c(s)}M \to T_{c(t)}M\) for parallel
transport along \(c\) from \(s\) to \(t\). Its inverse
\(\parallel_{t}^{0}(c)\colon T_{c(t)}M \to T_{c(0)}M\) is the
\emph{backward parallel transport}. The covariant derivative of a
vector field \(Y\) along \(c\) can be computed via
\begin{equation}\label{eq:backward}
  \frac{\nabla}{\mathrm{d}t}\bigg|_{t=0} Y(t)
  = \frac{\mathrm{d}}{\mathrm{d}t}\bigg|_{t=0}
  \bigl(\parallel_{t}^{0}(c) Y(t)\bigr)
\end{equation}
where the right-hand side is the ordinary derivative of the curve
\(\mathbb{R} \to T_{c(0)}M,\ t \mapsto \parallel_{t}^{0}(c) Y(t)\).

\begin{definition}[{cf.~\cite[Ch.~3]{JW1}}]
Let \(U \subset V\) be a star-shaped open neighbourhood of \(0\) in a
vector space \(V\), equipped with a Riemannian metric written in
geodesic normal coordinates at \(0\). The \emph{backward parallel
transport map} \(\Phi^{-1}\colon U \to \GL(V)\) assigns to \(X \in U\)
the backward parallel transport along the ray \(\gamma_X(t) := tX\):
\[
  \Phi^{-1}(X) := \parallel_{1}^{0}(\gamma_X)\colon T_XU \to T_0U .
\]
Using the canonical identifications \(T_XU \cong T_0U \cong V\), we
regard \(\Phi^{-1}(X)\) as an element of \(\GL(V)\).
\end{definition}

Then \(X \mapsto \Phi^{-1}(X)\) is smooth and, when \(V := T_pM\) with
the pulled-back metric \(\exp_p^{*}g\), it is the inverse of the
forward transport \(\Phi\) from \eqref{eq:Phi}. Moreover, if we view
\(Y \in V\) as the constant vector field \(Y_X=(X,Y)\) on \(U\) via the
trivialisation \(TU \cong U \times V\), \eqref{eq:backward} yields the
asymptotic expansion
\begin{equation}\label{%
eq:Taylor_expansion_of_the_backward_parallel_transport}
  \Phi^{-1}(X) Y \underset{X \to 0}{\sim}
  \sum_{k=0}^{\infty} \frac{1}{k!}\,\nabla^{k}_{X,\ldots,X} Y\big|_{0}.
\end{equation}

To describe the Taylor coefficients by noncommutative polynomials,
define \(\tilde Q_k \in
\scrA_{\text{univ}}(\scrR^0,\scrR^1,\scrR^2,\ldots)\) recursively by
\(\tilde Q_0=\Id_V\), \(\tilde Q_1=0\), and for \(k \ge 0\),
\begin{equation}\label{eq:bausteine}
  \tilde Q_{k+2} := -\frac{k+1}{k+3}
  \sum_{j=0}^{k} \binom{k}{j}\, \scrR^{j}\,\tilde Q_{k-j}.
\end{equation}

\begin{proposition}\label{%
p:taylor_series_of_the_backward_parallel_transport}
Let \(U \subset V\) be as above. Then the polynomials \(\tilde Q_k\)
satisfy
\begin{equation}\label{eq:formel_fuer_die_nabla_ks}
  \nabla^{k}_{X,\ldots,X} Y\big|_{0}
  = \tilde Q_k\bigl(\scrR^0(X),\scrR^1(X),\ldots\bigr)\,Y\big|_{0}.
\end{equation}
Equivalently, \(\frac{1}{k!}\,\tilde Q_k\bigl(\scrR^0|_0,\scrR^1|_0,
\ldots\bigr)\) is the \(k\)th coefficient of the Taylor expansion of
\(\Phi^{-1}\) in
\eqref{eq:Taylor_expansion_of_the_backward_parallel_transport}.
\end{proposition}

\begin{proof}
For \(k=0\) the claim is clear. For \(k\ge1\), fix \(X,Y\in V\) and set
\(\gamma(t):=tX\). Let \(J_Y\) be the unique Jacobi field along
\(\gamma\) with initial data \(J_Y(0)=0\) and \(J_Y^{1)}(0)=Y\).
Then
\begin{equation}\label{eq:Jacobi_equation}
  J_Y^{2)} \equiv -\scrR_\gamma\,J_Y ,
\end{equation}
where \(J_Y^{m)}\) denotes the \(m\)th covariant \(t\)-derivative and
\(\scrR_\gamma\) the Jacobi operator. Moreover, for \(k\ge0\),
\begin{align}
\label{eq:k-ter_Jacobi_on_a_geodesic}
  \scrR_{\gamma}^{k)}\big|_{t=0} &= \scrR^{k}(X)\big|_{0},\\
\label{eq:a_special_case_of_the_product_rule}
  J_Y^{k+1)}\big|_{t=0} &= (k+1)\,\nabla^{k}_{X,\ldots,X}Y\big|_{0}.
\end{align}
From these, for \(k\ge2\),
\[
  \begin{aligned}
  \nabla^{k}_{X,\ldots,X}Y\big|_{0}
    &= \frac{1}{k+1} J_Y^{k+1)}\big|_{t=0}
     = -\frac{1}{k+1} (\scrR_\gamma J_Y)^{k-1)}\big|_{t=0} \\
    &= -\frac{1}{k+1}\sum_{j=0}^{k-1} \binom{k-1}{j}\,
        \scrR_{\gamma}^{j)}\,J_Y^{k-1-j)}\big|_{t=0} \\
    &= -\frac{1}{k+1}\sum_{j=0}^{k-1} (k-1-j)\binom{k-1}{j}\,
        \scrR^{j}(X)\,\nabla^{k-2-j}_{X,\ldots,X}Y\big|_{0} \\
    &= -\frac{k-1}{k+1}\sum_{j=0}^{k-2} \binom{k-2}{j}\,
        \scrR^{j}(X)\,\nabla^{k-2-j}_{X,\ldots,X}Y\big|_{0},
  \end{aligned}
\]
which matches \eqref{eq:bausteine} together with the induction
hypothesis for \eqref{eq:formel_fuer_die_nabla_ks}.
\end{proof}

For example,
\[
  \tilde Q_2 = -\tfrac{1}{3}\scrR^0,\qquad
  \tilde Q_3 = -\tfrac{1}{2}\scrR^1,\qquad
  \tilde Q_4 = -\tfrac{3}{5}\scrR^2 + \tfrac{1}{5}\scrR^0\,\scrR^0 .
\]
Hence the Taylor polynomial of order four of the backward parallel
transport is (cf.~\cite[p.~332]{Gray})
\begin{equation}\label{eq:Taylor_polynomial_of_Phi-1}
  \begin{aligned}
  \Phi^{-1}(X)Y
    &= Y - \tfrac{1}{6}\scrR^0(X)Y - \tfrac{1}{12}\scrR^1(X)Y \\
    &\quad - \tfrac{1}{40}\scrR^2(X)Y
           + \tfrac{1}{120}\scrR^0(X)\,\scrR^0(X)Y .
  \end{aligned}
\end{equation}

To obtain a nonrecursive description of the \(\tilde Q_k\), set
\begin{equation}\label{eq:scrR_rescaled}
  \bar{\scrR}^{\,j} := -\frac{1}{j!}\,\scrR^{j}
\end{equation}
and note the canonical algebra anti-involution
\(*\colon \scrA_{\text{univ}}(\scrR^0,\scrR^1,\scrR^2,\ldots) \to
        \scrA_{\text{univ}}(\scrR^0,\scrR^1,\scrR^2,\ldots)\)
defined by \(\scrR^{i*}=\scrR^{i}\) and \((PQ)^* = Q^*P^*\). For a
sequence \(I=(i_1,\ldots,i_r)\) of nonnegative integers, write
\(\bar{\scrR}^{\,I}:=\bar{\scrR}^{\,i_1}\cdots\bar{\scrR}^{\,i_r}\) and
\begin{equation}\label{eq:Pi_I}
  \Pi_I := (i_1+2)(i_1+3)(i_1+i_2+4)(i_1+i_2+5)\cdots
           (i_1+\cdots+i_r+2r)(i_1+\cdots+i_r+2r+1).
\end{equation}
Then, either by \eqref{eq:bausteine} or directly from
\cite[Lemma~3.1]{JW1},
\begin{equation}\label{eq:explicit_formula_for_Q_k_tilde}
  \tilde Q_k = \sum_{\deg(I)=k}
  \frac{k!}{\Pi_I}\,\bar{\scrR}^{\,I*}.
\end{equation}

\subsection{Proof of Theorem~\ref{th:taylor_series}}
\label{se:taylor_series}

Because the Levi–Civita connection is metric (\(\nabla g=0\)),
\[
  g_X(Y,Z)=\langle \Phi^{-1}(X)Y, \Phi^{-1}(X)Z\rangle
          =\langle \Phi^{-1}(X)^{*}\,\Phi^{-1}(X)Y, Z\rangle
\]
where \(\langle\cdot,\cdot\rangle:=g_0\) and \(^*\) denotes the adjoint
with respect to \(\langle\cdot,\cdot\rangle\). Define, for \(k\ge0\),
\begin{equation}\label{eq:def_of_Q}
  Q_k:=\sum_{j=0}^{k}\binom{k}{j}\,\tilde Q_j^{*}\,\tilde Q_{k-j}
\end{equation}
where \(\tilde Q_k\) are given recursively by~\eqref{eq:bausteine} or
explicitly by~\eqref{eq:explicit_formula_for_Q_k_tilde}, and \(^*\) is
the canonical algebra anti-involution on
\(\scrA_{\text{univ}}(\scrR^0,\scrR^1,\scrR^2,\ldots)\) characterized by
\((PQ)^*=Q^*P^*\) and \((\scrR^i)^*=\scrR^i\). By
Proposition~\ref{p:taylor_series_of_the_backward_parallel_transport} and
the Cauchy product for Taylor series, \eqref{eq:def_of_Q} yields the
coefficients
\[
  \frac{1}{k!} Q_k\bigl(\scrR^0(X),\scrR^1(X),\ldots\bigr)
\]
in the Taylor expansion of the metric tensor in geodesic normal
coordinates stated in Theorem~\ref{th:taylor_series}.\qed

For example,
\[
  \begin{aligned}
    Q_2&=-\tfrac{2}{3}\,\scrR^0, &\quad
    Q_3&=-\,\scrR^1, \\
    Q_4&=-\tfrac{6}{5}\,\scrR^2+\tfrac{16}{15}\,\scrR^0\,\scrR^0, &\quad
    Q_5&=-\tfrac{4}{3}\,\scrR^3+\tfrac{8}{3}\bigl(\scrR^1\scrR^0+\scrR^0\scrR^1\bigr)
  \end{aligned}
\]
which gives the Taylor expansion \eqref{eq:Taylor_of_deg_5} of the
metric tensor (cf.~\cite[p.~336]{Gray}). Using
\eqref{eq:explicit_formula_for_Q_k_tilde}, we also have
\begin{equation}\label{eq:explicit_formula_for_Q_k}
  Q_k=\sum_{j=0}^{k}\sum_{\substack{\deg(I)=j\\ \deg(J)=k-j}}
  \frac{k!}{\Pi_I\,\Pi_J}\ \bar{\scrR}^{\,J}\,\bar{\scrR}^{\,I*}
\end{equation}
where the rescaled variables \(\bar{\scrR}^{\,j}\) are defined in
\eqref{eq:scrR_rescaled}.

\bigskip
\begin{corollary}
\begin{equation}\label{eq:Q_k}
  Q_{k+2}=c_k\,\scrR^{k}+\text{terms involving only }
  \scrR^{0},\ldots,\scrR^{k-1}
\end{equation}
with \(c_k=-2\,\frac{k+1}{k+3}\).
\end{corollary}

\noindent\emph{Sketch of proof.}
In \eqref{eq:def_of_Q} the leading term in \(\scrR^k\) comes from
\(j=0\) and \(j=k\): \(\tilde Q_0^*\,\tilde Q_k+\tilde Q_k^*\,\tilde Q_0
=2\,\tilde Q_k\) since \(\tilde Q_0=\Id\) and \(\tilde Q_k^*=\tilde
Q_k\) at top degree. The recursion \eqref{eq:bausteine} gives
\(\tilde Q_{k+2}\) leading term \(-\frac{k+1}{k+3}\scrR^{k}\);
multiplying by \(2\) yields \(c_k\) as stated.

\medskip
The Taylor expansion of the metric in geodesic normal coordinates is
also proved in~\cite{Gray} by a similar method; that approach does not
invoke the Jacobi equation. See also~\cite{MSV} for another derivation.

\section{Weyl's construction of irreducible representations of the
general linear group}
\label{se:YS}

Following Fulton–Harris \cite{Ful,FulH} and \cite[Ch.~4]{FKWC}, we
briefly review Young diagrams and tableaux, the associated symmetrizers
and projectors on tensor spaces, and their relation to irreducible
representations of the general linear group via Schur functors.

A partition \(\lambda_1 \ge \cdots \ge \lambda_k > 0\) of an integer
\(d\) can be depicted as a Young frame: an arrangement of \(d\) boxes
aligned from the left in \(k\) rows of lengths \(\lambda_i\) (top to
bottom). For example, the frame corresponding to \((5,3,2)\) is
\[
 \begin{array}{|c|c|c|c|c|c|c|}
  \cline{1-5}
   & & & & \\
  \cline{1-5}
   & & & \multicolumn{2}{c}{\;\;}\\
  \cline{1-3}
   & & \multicolumn{3}{c}{\;\;}\\
  \cline{1-2}
 \end{array}
\]
Filling the boxes with \(d\) distinct numbers \(n_1,\dots,n_d\) yields a
\emph{Young tableau} of shape \(\lambda\) (cf.~\cite{Ful}). For example,
\begin{equation}\label{eq:normal_standard_Young_tableau}
 T=\begin{array}{|c|c|c|c|c|}
     \cline{1-5}
     1 & 10 & 9 & 2 & 5\\
     \cline{1-5}
     8 & 7 & 4 & \multicolumn{1}{c}{\;\;} \\
     \cline{1-3}
     3 & 6 & \multicolumn{2}{c}{\;\;}\\
     \cline{1-2}
   \end{array}
\end{equation}
is a tableau of shape \((5,3,2)\). For simplicity we assume
\(\{n_1,\dots,n_d\}=\{1,\dots,d\}\). When these numbers appear
left-to-right in each row and top-to-bottom across rows, the tableau is
\emph{normal}: the entries \(1,\dots,\lambda_1\) occupy the first row,
\(\lambda_1+1,\dots,\lambda_1+\lambda_2\) the second, and so on.

Let \(V\) be a real vector space with dual \(V^*\). The symmetric group
\(S_d\) acts on the right by
\[
  X_1\otimes\cdots\otimes X_d\cdot\sigma
  := X_{\sigma(1)}\otimes\cdots\otimes X_{\sigma(d)}
\]
on \(\bigotimes^d V\), and hence on the left by
\[
  (\sigma\cdot\lambda)(X_1,\dots,X_d)
  := \lambda(X_{\sigma(1)},\dots,X_{\sigma(d)})
\]
on \(\bigotimes^d V^*\).

Fix a tableau \(T\) of shape \(\lambda\), and let \(S_r\) and \(S_c\) be
the subgroups of \(S_d\) preserving its rows and columns, respectively.
The row symmetrizer and column antisymmetrizer are
\begin{align}\label{eq:def_r_lambda}
  r_T \colon \bigotimes^d V^* \to \bigotimes^d V^*,\quad
  \lambda \mapsto \sum_{\sigma\in S_r} \lambda\cdot\sigma \\
\label{eq:def_c_lambda}
  c_T \colon \bigotimes^d V^* \to \bigotimes^d V^*,\quad
  \lambda \mapsto \sum_{\sigma\in S_c} (-1)^{|\sigma|}\,\lambda\cdot\sigma
\end{align}
and the associated Young symmetrizers on \(\bigotimes^d V^*\) are
\begin{equation}\label{eq:def_S_*_lambda}
  S_T := r_T \circ c_T,\qquad
  S_T^\star := c_T \circ r_T
\end{equation}
(cf.~\cite[p.~46, (4.2)]{FulH}). Their images,
\[
  \bbS_T V^* := S_T\bigl(\bigotimes^d V^*\bigr),\qquad
  \bbS_T^\star V^* := S_T^\star\bigl(\bigotimes^d V^*\bigr)
\]
are \(\GL(V)\)-modules. After complexifying,
\((\bbS_T V^*)_\C\) and \((\bbS_T^\star V^*)_\C\) are irreducible
polynomial \(\GL(V_\C)\)-modules with highest weight \(\lambda\). The
maps \(c_T\) and \(r_T\) give explicit \(\GL(V)\)-equivariant
isomorphisms \(\bbS_T V^* \cong \bbS_T^\star V^*\). The assignments
\(V \mapsto \bbS_T V^*\) and \(V \mapsto \bbS_T^\star V^*\) are the
(covariant) Schur functor and its dual associated with \(T\).

By Schur’s Lemma there exists a constant \(h_\lambda\in\N_{\ge1}\),
depending only on the frame, such that \(P_T:=\frac{1}{h_\lambda}S_T\)
and \(P_T^\star:=\frac{1}{h_\lambda}S_T^\star\) are projectors (the
Young projectors). For each box of the frame, its \emph{hook length} is
the number of boxes weakly to its right in the same row plus the number
weakly below it in the same column minus one; then \(h_\lambda\) is the
product of all hook lengths over the diagram.

Following \cite[Ch.~15.5]{FulH}, there is another characterization of
\(\bbS_T^\star V^*\). Let \(\mu_1 \ge \cdots \ge \mu_\ell\) be the
conjugate partition (column lengths). Then
\(\bbS_T^\star V^* \subset \Lambda^{\mu_1}V^* \otimes \cdots \otimes
\Lambda^{\mu_\ell}V^*\). Moreover, if \(T\) is the transpose of a normal
tableau, the numbers \(1,\dots,d\) are written (top to bottom, left to
right) into the boxes of \(T\): \(1,\dots,\mu_1\) fill the first column,
\(\mu_1+1,\dots,\mu_1+\mu_2\) the second, etc.

\begin{theorem}\label{th:char1}
We have
\begin{equation}\label{eq:char1}
  \bbS_T^\star V^* = \bigcap_{i<j} \Kern(\ell_{ij}^\star)
\end{equation}
where \(\ell_{ij}^\star\) is the dual of the canonical map
\(\Lambda^{\mu_i+1}V \otimes \Lambda^{\mu_j-1}V \to
  \Lambda^{\mu_i}V \otimes \Lambda^{\mu_j}V\):
\begin{equation}\label{eq:lij}
  \begin{aligned}
   v_1 \wedge \cdots \wedge v_{\mu_i+1}
   &\otimes v_{\mu_i+2} \wedge \cdots \wedge v_{\mu_i+\mu_j} \longmapsto\\
  & \sum_{a=1}^{\mu_i+1} (-1)^{a+\mu_i+1}
   v_1 \wedge \cdots \wedge \hat v_a \wedge \cdots \wedge v_{\mu_i+1}
   \otimes v_a \wedge v_{\mu_i+2} \wedge \cdots \wedge v_{\mu_i+\mu_j}
  \end{aligned}
\end{equation}
for \(i<j\).
\end{theorem}

For the shape \((k+2,2)\), the maps \(\ell_{1,2}^\star\) and
\(\ell_{1,3}^\star\) encode the first and second Bianchi identities; i.e.,
\[
  \bbS^\star_{\scaleto{\begin{array}{|c|c|c|c|c|c|}
   \cline{1-5}
   1 & 3 & 5 & \cdots & k+4\\
   \cline{1-5}
   2 & 4 & \multicolumn{2}{c}{\;\;} \\
   \cline{1-2}
  \end{array}}{15pt}} V^* = \scrC_k^\star(V)
\]
is the space of linear algebraic \(k\)-jets of the curvature tensor (see
Definition~\ref{de:algebraic_curvature_jets}).

A parallel description of \(\bbS_T V^*\) is known (see
\cite[Ch.~8.3, Ex.~10]{Ful}). Assume now that \(T\) is normal.

\begin{theorem}\label{th:char2}
We have
\begin{equation}\label{eq:char2}
  \bbS_T V^* = \bigcap_{i<j} \Kern(\ell_{ij})
\end{equation}
where \(\ell_{ij}\) is the dual of the canonical map
\(\Sym^{\lambda_i+1}V \otimes \Sym^{\lambda_j-1}V \to
  \Sym^{\lambda_i}V \otimes \Sym^{\lambda_j}V\):
\[
  \begin{aligned}
   v_1 \odot \cdots \odot v_{\lambda_i+1}
   &\otimes v_{\lambda_i+2} \odot \cdots \odot v_{\lambda_i+\lambda_j} \longmapsto\\
   & \sum_{a=1}^{\lambda_i+1}
   v_1 \odot \cdots \odot \hat v_a \odot \cdots \odot v_{\lambda_i+1} 
   \otimes v_a \odot v_{\lambda_i+2} \odot \cdots \odot v_{\lambda_i+\lambda_j}
  \end{aligned}
\]
\end{theorem}
Here \(\odot\) denotes the symmetric product. In particular,
\[
  \bbS_{\scaleto{\begin{array}{|c|c|c|c|c|c|}
   \cline{1-5}
   1 & 3 & 5 & \cdots & k+4\\
   \cline{1-5}
   2 & 4 & \multicolumn{2}{c}{\;\;} \\
   \cline{1-2}
  \end{array}}{15pt}} V^* = \scrC_k(V)
\]
(see ~\ref{eq:formal_Jacobi_operator}). For two rows, the
proof of Theorem~\ref{th:char2} follows directly from Weyl’s dimension
formula via a short exact sequence similar to \eqref{eq:ses};
similarly, Theorem~\ref{th:char1} follows from a single short exact
sequence when the diagram has two columns.

\section*{ACKNOWLEDGEMENTS}
I would like to thank Gregor Weingart for valuable remarks.

\bibliographystyle{amsplain}

\end{document}